\documentclass[11pt,a4paper]{amsart}
\usepackage{xcolor}
\usepackage[T1]{fontenc}
\usepackage{lmodern}
\usepackage{amsmath,amssymb,amsthm,mathtools,microtype}
\usepackage[margin=28mm]{geometry}
\usepackage{enumitem}
\usepackage[hidelinks]{hyperref}
\hypersetup{pdftitle={A weak Nullstellensatz for commuting skew variables over division rings},pdfauthor={},pdfsubject={Skew polynomial rings, simple modules, and twisted fixed division rings}}
\setlist[enumerate,1]{label=\textup{(\roman*)},leftmargin=2.3em,itemsep=3pt,topsep=5pt}
\allowdisplaybreaks[1]
\numberwithin{equation}{section}
\newtheorem{theorem}{Theorem}[section]
\newtheorem{proposition}[theorem]{Proposition}
\newtheorem{lemma}[theorem]{Lemma}
\newtheorem{corollary}[theorem]{Corollary}
\theoremstyle{definition}
\newtheorem{definition}[theorem]{Definition}
\newtheorem{example}[theorem]{Example}
\newtheorem{remark}[theorem]{Remark}

\newcommand{\N}{\mathbb N}
\newcommand{\HH}{\mathbb H}
\newcommand{\CC}{\mathbb C}
\newcommand{\RR}{\mathbb R}
\newcommand{\FF}{\mathbb F}
\newcommand{\ba}{\mathbf a}
\newcommand{\bb}{\mathbf b}
\newcommand{\bc}{\mathbf c}
\newcommand{\Hcond}{\mathcal H}
\newcommand{\Cen}{\mathrm C}
\DeclareMathOperator{\Aut}{Aut}
\DeclareMathOperator{\End}{End}
\DeclareMathOperator{\Int}{Int}
\DeclareMathOperator{\Ann}{Ann}

\DeclareMathOperator{\lm}{lm}
\DeclareMathOperator{\lc}{lc}
\DeclareMathOperator{\Span}{span}
\DeclareMathOperator{\id}{id}
\newcommand{\doi}[1]{\href{https://doi.org/#1}{\nolinkurl{#1}}}

\title[Hilbert's Nullstellensatz for  skew polynomials over division rings]{Weak Hilbert's Nullstellensatz for  skew polynomials over division rings}

\date{}
\author{ Dinh Van Hoang}
\address[Dinh Van Hoang]{Faculty of Advanced Education, Ho Chi Minh City University of Technology and Engineering, Vietnam }
\email{hoangdv@hcmute.edu.vn}

\author{Vo Ngoc Thieu}
\address[Vo Ngoc Thieu]{Department of Computer Science, University of Bath, United Kingdom}
\email{ntv22@bath.ac.uk}

\author{Phan Thanh Toan}
\address[Phan Thanh Toan]{Analytical and Algebraic Methods in Optimization Research Group, Faculty of Mathematics and Statistics, Ton Duc Thang University, Ho Chi Minh City, Vietnam}
\email{phanthanhtoan@tdtu.edu.vn}

\begin{document}
\begin{abstract}
We establish a weak Nullstellensatz for skew polynomial rings
$R_n=D[x_1,\ldots,x_n;\sigma]$, where $D$ is a division ring and
the commuting variables induce the same automorphism $\sigma$ of $D$.
For $n\ge2$, we characterize the pairs $(D,\sigma)$ for which every
maximal left ideal is generated by $x_1-a_1,\ldots,x_n-a_n$, with
$\sigma(a_j)a_i=\sigma(a_i)a_j$. The criterion consists of right
$\sigma$-algebraic closure of $D$ and right algebraic closure of the
centralizers of all commuting $(n-2)$-tuples in each twisted fixed
division ring $E_a=\{c\in D:\sigma(c)a=ac\}$ for $0\ne a\in D$.
\end{abstract}
\maketitle
%\tableofcontents
\section{Introduction and main results}\label{sec:intro}

The weak Nullstellensatz identifies maximal ideals in a polynomial ring
over an algebraically closed field with points of affine space. Over a
division ring, this correspondence must be formulated using one-sided
ideals, and even the choice of evaluation points depends on the
commutation relations between the variables and the coefficients. This
paper treats commuting skew variables that induce the same automorphism of
the coefficient division ring.

For central variables, the theorem of Amitsur and Small
\cite{AmitsurSmall} establishes finite-dimensionality of simple modules
over $D[x_1,\ldots,x_n]$. Alon and Paran \cite{AlonParan} proved the
central quaternionic Nullstellensatz. Aryapoor \cite{AryCentral, Aryapoor2024}
subsequently characterized the division rings satisfying the central
Nullstellensatz by algebraic closure conditions on centralizers of
commuting tuples and gave the explicit Nullstellensatz for the division ring of quaternions. The distinction between nonzero contraction and
maximal contraction of left ideals is important in this subject; see
Chapman and Paran \cite{ChapmanParan}.

The skew setting introduces a further difficulty. One-variable division
and right evaluation are governed by Ore's theory \cite{Ore} and its
developments, including \cite{LamLeroy}. For the several-variable ring
considered here, evaluation at compatible tuples also appears in the
work of Alon, Behajaina, and Paran \cite{AlonBehajainaParan}.
Aryapoor's skew Amitsur--Small theorem
\cite[Theorem~1.3 of the preprint version]{ArySkew} gives
finite-dimensionality under an additional hypothesis that the
one-variable skew polynomial ring contains a nonconstant central
polynomial. We do not impose that hypothesis. Instead, we first analyze
finite-dimensional simple modules and then show that the same criterion
controls all simple modules. The Puiseux-series example in
Section~\ref{subsec:puiseux} has no nonconstant central polynomial in
its one-variable skew polynomial ring.

\subsection{Notation and closure conditions}
All rings are associative with identity, and all modules are unital left
modules unless explicitly stated otherwise. Let $D$ be a division ring
and $\sigma\in\Aut(D)$. We write
\[
R_n=R_n(D,\sigma)=D[x_1,\ldots,x_n;\sigma]
\]
for the free left $D$-module with basis $x^\alpha$, $\alpha\in\N^n$,
where $\N=\{0,1,\ldots\}$, and multiplication
\begin{equation}\label{eq:mult}
(d x^\alpha)(e x^\beta)
   =d\,\sigma^{|\alpha|}(e)x^{\alpha+\beta}.
\end{equation}
Thus $x_i x_j=x_j x_i$ and $x_i d=\sigma(d)x_i$.
For a tuple $\bb$ in a division ring $E$, let
\[
\Cen_E(\bb)=\{c\in E:cb_j=b_jc\text{ for every }j\}.
\]
The centralizer of the empty tuple is understood to be $E$.

\begin{definition}\label{def:closure}
	\begin{enumerate}
		\item A division ring $E$ is \emph{right algebraically closed} if every
		positive-degree polynomial $p(z)=\sum_j e_jz^j\in E[z]$, with $z$
		central, has a right root: $\sum_j e_jb^j=0$ for some $b\in E$.
		\item Following \cite{AryCentral}, $E$ is \emph{centrally algebraically closed}
		if $\Cen_E(\bb)$ is right algebraically closed for every finite
		pairwise commuting tuple $\bb$ in $E$, including the empty tuple.
		\item The pair $(D,\sigma)$ is \emph{right $\sigma$-algebraically closed}
		if every positive-degree polynomial in $D[t;\sigma]$ has a right root.
		Here a right root $a$ of $p$ means that $p\in D[t;\sigma](t-a)$.
	\end{enumerate}

\end{definition}

For central polynomials, right algebraic closure is equivalent to left
algebraic closure: both are equivalent to the assertion that all
irreducible polynomials have degree one. These notions do not require
the coefficient division ring to be commutative.

For $a\in D^\times$, define the \emph{twisted fixed division ring}
\begin{equation}\label{eq:Ea}
E_a=E_a(D,\sigma)
 =\{c\in D:\sigma(c)a=ac\}
 =D^{\theta_a},
\qquad \theta_a=\Int(a^{-1})\circ\sigma,
\end{equation}
where $\Int(u)(d)=udu^{-1}$. The admissible points are
\begin{equation}\label{eq:admissible}
D^{n,\sigma}
 =\{\ba\in D^n:\sigma(a_j)a_i=\sigma(a_i)a_j\text{ for all }i,j\}.
\end{equation}
For $\ba\in D^{n,\sigma}$, the operators
\begin{equation}\label{eq:point-action}
T_i^{\ba}(d)=\sigma(d)a_i
\end{equation}
define an $R_n$-module on the left $D$-space $D$. We denote this module
by $D_{\ba}$ and define
\[
f(\ba)=f(T^{\ba})1,
\qquad
\mathfrak m_{\ba}=\sum_{i=1}^nR_n(x_i-a_i).
\]
Section~\ref{sec:prelim} proves that $\mathfrak m_{\ba}$ is precisely
the left ideal of polynomials vanishing at $\ba$. For $S\subseteq R_n$,
let
\[
V_\sigma(S)=\{\ba\in D^{n,\sigma}:f(\ba)=0\text{ for all }f\in S\}.
\]

\noindent For $X\subseteq D^{n,\sigma}$, define
\[
I_\sigma(X)=\{f\in R_n:f(\ba)=0\text{ for every }\ba\in X\}.
\]
By Proposition~\ref{prop:pointideal}, this is the left ideal
$\bigcap_{\ba\in X}\mathfrak m_{\ba}$.

\medskip
For $n=1$, let $\Hcond_1(D,\sigma)$ mean right
$\sigma$-algebraic closure. For $n\ge2$, let
$\Hcond_n(D,\sigma)$ denote the following two conditions:
\begin{itemize}
	\item[(1)] $(D,\sigma)$  is right $\sigma$-algebraically closed;
	\item[(2)] $\Cen_{E_a}(b_1,\ldots,b_{n-2})$  is right algebraically closed {for every }$a\in D^\times${ and every pairwise commuting}
	 $(b_1,\ldots,b_{n-2})\in E_a^{n-2}$.
\end{itemize}

Padding tuples with zeros shows that $\Hcond_n$ implies $\Hcond_m$
whenever $1\le m\le n$.

\subsection{The main results}

\begin{theorem}[Weak Nullstellensatz]\label{thm:main}
Let $D$ be a division ring, let $\sigma\in\Aut(D)$, and let $n\ge1$.
The following conditions are equivalent.
\begin{enumerate}
\item\label{main:H} The condition $\Hcond_n(D,\sigma)$ holds.
\item\label{main:all} Every simple left $R_n$-module has left
$D$-dimension one.
\item\label{main:finite} Every finite-dimensional simple left
$R_n$-module has left $D$-dimension one.
\item\label{main:ideals} The maximal left ideals of $R_n$ are exactly
$\mathfrak m_{\ba}$, where $\ba\in D^{n,\sigma}$.
\item\label{main:zeros} Every proper left ideal $I$ of $R_n$ satisfies
$V_\sigma(I)\ne\varnothing$.
\item\label{main:certificate} For every finite family
$f_1,\ldots,f_r\in R_n$,
\begin{equation}\label{eq:certificate}
V_\sigma(\{f_1,\ldots,f_r\})=\varnothing
\quad\Longleftrightarrow\quad
\sum_{j=1}^r g_jf_j=1
\quad\text{for some }g_1,\ldots,g_r\in R_n.
\end{equation}
\end{enumerate}
Moreover, if $\Hcond_n(D,\sigma)$ holds, then for every left ideal $I\subseteq R_n$,
\begin{equation}\label{eq:radical}
	I_\sigma\bigl(V_\sigma(I)\bigr)
	=\bigcap_{\substack{L\supseteq I\\L\text{ maximal left}}}L.
\end{equation}
%No assumption is made on the dimension of $D$ over its center or, in advance, on the dimension of arbitrary simple modules.
\end{theorem}

The preceding criterion leads to an existence question: can an arbitrary
pair $(D,\sigma)$ be embedded into one satisfying $\Hcond_n$ for every
$n$?  An embedding of pairs is an injective unital homomorphism
$\iota:D\to E$ together with an automorphism
$\widetilde{\sigma}\in\Aut(E)$ such that
$\widetilde{\sigma}\circ\iota=\iota\circ\sigma$.
After identifying $D$ with its image, this means that
$\widetilde{\sigma}|_D=\sigma$. We answer this question by combining inner twisting with the embedding theorem of Aryapoor
\cite[Theorem~5.2]{AryCentral}.
\begin{theorem}[Equivariant embedding]\label{thm:embedding}
	Let $D$ be a division ring and let $\sigma\in\Aut(D)$. There exist
	a division ring $E\supseteq D$ and an element $u\in E^\times$ such that,
	with $\widetilde{\sigma}=\Int(u)$, the following properties hold:
	\begin{enumerate}
		\item $E$ is centrally algebraically closed.
		\item $\widetilde{\sigma}|_D=\sigma$, and $\widetilde{\sigma}$ has
		the same order as $\sigma$, whether finite or infinite.
		\item $(E,\widetilde{\sigma})$ is right
		$\widetilde{\sigma}$-algebraically closed.
		\item For every $a\in E^\times$,
		\begin{equation}\label{eq:embedding-fixed}
			E_a(E,\widetilde{\sigma})=\Cen_E(u^{-1}a),
		\end{equation}
		and this division ring is centrally algebraically closed.
	\end{enumerate}
	More precisely, the construction gives a chain
	$D\subseteq K\subseteq E$, where $K=D(u;\sigma)$ is the division
	ring of fractions of $D[u;\sigma]$ and $Z(E)=Z(K)$.
\end{theorem}

\begin{corollary}\label{cor:embedding-alln}
	Every pair $(D,\sigma)$ has an extension
	$(E,\widetilde{\sigma})$ such that
	$\Hcond_n(E,\widetilde{\sigma})$ holds for every $n\ge1$.
	The automorphism $\widetilde{\sigma}$ may be chosen inner and of the
	same order as $\sigma$. In particular, for every $n\ge1$, all six
	conditions of Theorem~\ref{thm:main} hold for
	$R_n(E,\widetilde{\sigma})$.
\end{corollary}

For $n=2$, the centralizer condition (2) in $\mathcal{H}_2(D,\sigma)$  simply requires
every $E_a$ to be right algebraically closed. For all numbers of
variables, it takes the following form.

\begin{corollary}\label{cor:alln}
	The  weak Nullstellensatz holds for $R_n(D,\sigma)$ for every $n\ge1$
	if and only if $(D,\sigma)$ is right $\sigma$-algebraically closed and
	$E_a$ is centrally algebraically closed for every $a\in D^\times$.
\end{corollary}
When $D$ is a field, $E_a=D^\sigma$ for all $a\ne0$, and the criterion
is independent of $n\ge2$. When $\sigma=\id$, it reduces to the
centralizer criterion for central polynomial rings. In general the
rings $E_a$ need not be commutative, and their centralizers cannot be
replaced by the fixed subfield of the center of $D$.

\smallskip
The conditions $\Hcond_n(D,\sigma)$ involve centralizers of commuting
tuples of increasing length as $n$ grows. It is therefore natural to ask
when the two-variable condition $\Hcond_2(D,\sigma)$ already implies
the conditions for every number of variables. The next theorem gives
two sufficient hypotheses on $D$ and $\sigma$.
\begin{theorem}[Centralizer reduction]\label{thm:compression}
	\begin{enumerate}
		\item If $D$ is centrally finite, then, for every $\sigma\in\Aut(D)$,
		\[
		\Hcond_2(D,\sigma)\ \Longrightarrow\
		\Hcond_n(D,\sigma)\quad\text{for all }n\geq1.
		\]
		\item Suppose $D$ is algebraic over its center $F$, the field $F$ is infinite and
		perfect, and some positive power of $\sigma$ is inner. 
		Then
		\[
		\Hcond_2(D,\sigma)\ \Longrightarrow\
		\Hcond_n(D,\sigma)\quad\text{for all }n\geq1.
		\]
	\end{enumerate}
	
\end{theorem}

\subsection{Outline of the proof}
The proof of Theorem~\ref{thm:main} begins with finite-dimensional simple
modules. After choosing a nonzero variable action $T_1$, the action
$T_iT_1^{-1}$ are commuting $D$-linear module endomorphisms. Restriction
to the one-variable skew polynomial ring and passage to an eigenspace
reduce the problem to a central polynomial ring over a twisted fixed
division ring $E_a$. The centralizer condition in $\Hcond_n$ then forces
the original module to have $D$-dimension one.

To treat an arbitrary simple module, we introduce an auxiliary central
variable $z$ and prove that every maximal left ideal of $R_n[z]$ meets
$D[z]$ nontrivially. A leading-coefficient argument, using suitable
polynomials over the prime field, establishes this contraction result.
It follows that, for every $E\in\End_{R_n}(M)$ and $v\in M$, the space
$D[E]v$ is finite-dimensional. The simple $D[z]$-module types occurring
in $M$ form a finite orbit under a nonzero variable action. This orbit
produces a finite-dimensional auxiliary module in two skew variables.
The finite-dimensional criterion forces the orbit to have length one
and its simple constituents to have $D$-dimension one. Applying this
splitting to $E=T_nT_1^{-1}$ allows descent to a centralizer with one
fewer skew variable. The closure conditions pass to the descended pair,
so induction proves sufficiency for arbitrary simple modules. Conversely,
finite-dimensional test modules built from irreducible polynomials over
the relevant centralizers prove necessity. The equivalences involving
maximal left ideals, common zeros, and finite certificates then follow
from the description of point modules.

The equivariant embedding theorem is proved by adjoining a skew
indeterminate $u$, passing to $D(u;\sigma)$, and applying the
center-preserving embedding theorem for centrally algebraically closed
division rings. Conjugation by $u$ extends $\sigma$, and preservation of
the center ensures preservation of its order. For centralizer reduction,
the centrally finite case follows from a dimension bound for fixed
division rings and the classification of centrally finite right
algebraically closed division rings. The algebraic case uses a fixed
implementing element for an inner power of $\sigma$ and a weighted
primitive-element argument to realize each required centralizer as
another twisted fixed division ring.

Section~\ref{sec:prelim} establishes the evaluation and module-theoretic
preliminaries. Section~\ref{sec:contraction} proves central-variable
contraction and local finite-dimensionality for endomorphisms, and
Section~\ref{sec:finite} treats finite-dimensional simple modules.
Section~\ref{sec:descent} develops finite-orbit splitting and centralizer
descent. Section~\ref{sec:proof} proves sufficiency for arbitrary simple
modules and necessity from finite-dimensional ones.
Section~\ref{sec:compression-tools} supplies the lemmas used in the
centralizer reduction theorem. Section~\ref{sec:mainProofs} completes the proofs
of all three main theorems and their stated corollaries. Finally,
Section~\ref{sec:applications} presents the field, central, Puiseux-series,
Frobenius, and quaternionic specializations and examples.

\section{Evaluation and preliminary module lemmas}\label{sec:prelim}

\subsection{Point modules and right evaluation}
The action of $c\in D$ on $D_{\ba}$ is left multiplication by $c$.
For $\ba\in D^{n,\sigma}$, the maps in \eqref{eq:point-action} satisfy
\[
T_i^{\ba}(cd)=\sigma(c)T_i^{\ba}(d),
\qquad
T_i^{\ba}T_j^{\ba}(d)
   =\sigma^2(d)\sigma(a_j)a_i.
\]
The admissibility equations therefore express exactly the commutativity
of the variable actions. This verifies all defining relations of $R_n$.

We define an action of \(R_n\) on the additive group of \(D\) by specifying the action of its generators:
\[ \begin{cases} c\cdot d=cd,&c\in D,\\ x_i\cdot d=T_i^{\mathbf a}(d)=\sigma(d)a_i. \end{cases} \]
To obtain a well-defined left \(R_n\)-module, the operators must respect the defining relations of \(R_n\). For example, the relation \(x_ic=\sigma(c)x_i\) requires
\[ x_i\cdot(cd)=\sigma(c)(x_i\cdot d), \]
which is exactly semilinearity. 
Likewise, \(x_ix_j=x_jx_i\) requires
\[ T_i^{\mathbf a}T_j^{\mathbf a} =T_j^{\mathbf a}T_i^{\mathbf a}. \]
Both conditions have been verified. Thus the action extends to all of \(R_n\), giving a left $R_n$ module structure of $D$ which we denote by 
\[ {D_{\mathbf a}}. \]
Moreover, $D_{\ba}$ is simple because any nonzero submodule is a nonzero left
$D$-subspace of the one-dimensional space $D$.

Define
\[
N_0(\ba)=1,
\qquad
N_{\alpha+e_i}(\ba)=\sigma(N_\alpha(\ba))a_i.
\]
Directly from the definition, for $\ba=(a,b)$,
\[
\begin{aligned}
	N_{1,0}&=a,
	&N_{0,1}&=b,\\
	N_{2,0}&=\sigma(a)a,
	&N_{0,2}&=\sigma(b)b,\\
	N_{1,1}&=\sigma(a)b=\sigma(b)a.
\end{aligned}
\]
Commutativity of the operators makes this recursion independent of the
chosen path through $\N^n$. Consequently,
\begin{equation}\label{eq:evaluation}
f(\ba)=\sum_\alpha d_\alpha N_\alpha(\ba)
\quad\text{if}\quad f=\sum_\alpha d_\alpha x^\alpha.
\end{equation}
In one variable, $N_0(a)=1$ and
\[
N_j(a)=\sigma^{j-1}(a)\cdots\sigma(a)a\qquad(j\ge1).
\]

\begin{proposition}\label{prop:pointideal}
For every $\ba\in D^{n,\sigma}$, the evaluation map
\[
\epsilon_{\ba}:R_n\longrightarrow D_{\ba}\equiv D,
\qquad f\longmapsto f\cdot1,
\]
is a surjective left $R_n$-module homomorphism with kernel
$\mathfrak m_{\ba}$. Distinct admissible points give distinct point
ideals. For $f,g\in R_n$,
\begin{equation}\label{eq:product}
(gf)(\ba)=g(T^{\ba})\bigl(f(\ba)\bigr).
\end{equation}
In particular, a left ideal certificate of $1$:
$$g_1f_1+\dots+g_rf_r=1$$precludes a common
admissible zero of $f_1,\dots,f_r$.
\end{proposition}
\begin{proof}
We first regard $R_n$ as its left regular module. For $h,f\in R_n$,
associativity of the action on $D_{\ba}$ gives
\[
\epsilon_{\ba}(hf)=(hf)\cdot1=h\cdot(f\cdot1)
                  =h\cdot\epsilon_{\ba}(f).
\]
Thus $\epsilon_{\ba}$ is a module homomorphism, even though evaluation
need not be a ring homomorphism. It is surjective because
$\epsilon_{\ba}(d)=d$ for every constant $d\in D$. Also,
$(x_i-a_i)\cdot1=0$ for each $i$. Hence
\[
\mathfrak m_{\ba}\subseteq\ker\epsilon_{\ba}.
\]

To prove the reverse inclusion, we show that every polynomial is
congruent to a constant modulo $\mathfrak m_{\ba}$. This is an
induction on the total degree of a monomial. The assertion is immediate
in degree zero. For a monomial $x^\alpha$ of positive degree, choose
$i$ with $\alpha_i>0$ and put $\beta=\alpha-e_i$. The defining
relations give
\[
x^\alpha=x^\beta(x_i-a_i)+x^\beta a_i
        =x^\beta(x_i-a_i)+\sigma^{|\beta|}(a_i)x^\beta.
\]
The first term belongs to the left ideal $\mathfrak m_{\ba}$, and
the second has smaller total degree. The induction hypothesis therefore
applies. Taking left $D$-linear combinations proves the assertion for
an arbitrary $f\in R_n$.

Write $f-c\in\mathfrak m_{\ba}$ with $c\in D$. Applying
$\epsilon_{\ba}$ gives $f(\ba)=c$, because
$\mathfrak m_{\ba}\subseteq\ker\epsilon_{\ba}$. In particular,
\[
f-f(\ba)\in\mathfrak m_{\ba}.
\]
If $f$ is in the kernel, this identity places $f$ in
$\mathfrak m_{\ba}$. Hence the kernel is exactly
$\mathfrak m_{\ba}$, and the induced isomorphism
$R_n/\mathfrak m_{\ba}\simeq D_{\ba}$ also shows that this left
ideal is maximal.

If $\mathfrak m_{\ba}=\mathfrak m_{\bb}$, then each $x_i-a_i$
belongs to the kernel of $\epsilon_{\bb}$. Its value at $\bb$ is
$b_i-a_i=0$, so $b_i=a_i$ for every $i$.

Finally, the action of a product is the composition of the two actions:
\[
(gf)(\ba)=(gf)\cdot1=g\cdot(f\cdot1)
                  =g(T^{\ba})\bigl(f(\ba)\bigr).
\]
If all $f_j$ vanish at $\ba$, every summand $g_jf_j$ also vanishes
there. Their sum cannot equal $1$, whose value is $1$. %This last argument uses only the module action and does not require an algebraic closure hypothesis.
\end{proof}

\begin{remark}\label{rem:nonadmissible}
Admissibility is essential, rather than a restriction added only for
convenience. For any tuple $\ba\in D^n$,
\begin{align*}
&(x_j-\sigma(a_j))(x_i-a_i)
 -(x_i-\sigma(a_i))(x_j-a_j)\\
&\hspace{35mm}=\sigma(a_j)a_i-\sigma(a_i)a_j.
\end{align*}
If $\ba$ is not admissible, the generated left ideal contains a nonzero
scalar and hence equals $R_n$. Also, evaluation is generally not a ring
homomorphism; \eqref{eq:product}, not multiplicativity in $D$, is the
identity used throughout.
\end{remark}

\begin{proposition}[Geometry and endomorphisms of point modules]\label{prop:geometry}
Suppose $a_1=a\ne0$ and write $a_i=ac_i$, with $c_1=1$.
The tuple $\ba=(a,ac_2,\dots,ac_n)$ is admissible if and only if $c_2,\ldots,c_n$ belong
to $E_a
 =\{c\in D:\sigma(c)a=ac\}
 $ and commute pairwise. In that case,
\begin{equation}\label{eq:point-end}
\End_{R_n}(D_{\ba})
 \simeq\Cen_{E_a}(c_2,\ldots,c_n)^{\mathrm{op}}.
\end{equation}
At the zero point the endomorphism ring is $D^{\mathrm{op}}$.
Furthermore, $D_{\ba}\simeq D_{\bb}$ if and only if there exists
$u\in D^\times$ such that
\begin{equation}\label{eq:conjugacy}
b_i=\sigma(u)a_i u^{-1}\qquad(1\le i\le n).
\end{equation}
\end{proposition}
\begin{proof}
We separate the assertions about points, endomorphisms, and
isomorphisms.

\smallskip\noindent
\emph{Admissibility.}
For $i\ge2$, the equation involving the first and the $i$th coordinates
is
\[
\sigma(ac_i)a=\sigma(a)ac_i.
\]
Since $\sigma(a)\ne0$, cancellation on the left reduces this to
$\sigma(c_i)a=ac_i$, which is precisely the condition $c_i\in E_a$.
Once these equations hold, the remaining compatibility equations become
\begin{align*}
\sigma(a_j)a_i
 &=\sigma(a)\sigma(c_j)ac_i
  =\sigma(a)a c_jc_i,\\
\sigma(a_i)a_j
 &=\sigma(a)\sigma(c_i)ac_j
  =\sigma(a)a c_ic_j.
\end{align*}
The common factor $\sigma(a)a$ is invertible. Thus equality holds
exactly when $c_ic_j=c_jc_i$. These computations are reversible,
which proves both directions of the admissibility assertion.

\smallskip\noindent
\emph{Endomorphisms.} Each $R_n$-endomorphism of $D_\ba$ is a left $D$-linear map from $D_\ba=D$ to itself. 
Each left $D$-linear map $F:D\to D$ is determined by $c=F(1)$, since
$F(d)=dF(1)=dc$. It is therefore the right multiplication map $r_c$.
For $d\in D$, the equations
\[
T_i^{\ba}r_c(d)=\sigma(d)\sigma(c)a_i,
\qquad
r_cT_i^{\ba}(d)=\sigma(d)a_ic
\]
show that $r_c$ commutes with the $i$th variable action if and only if
$\sigma(c)a_i=a_ic$. The first coordinate requires $c\in E_a$.
Using $a_i=ac_i$ and $\sigma(c)a=ac$, the other equations reduce to
\[
ac c_i=ac_ic,
\quad\text{or equivalently}\quad cc_i=c_ic.
\]
Hence the allowed coefficients are exactly
$\Cen_{E_a}(c_2,\ldots,c_n)$. The order of multiplication is reversed
under this correspondence: $r_c\circ r_d=r_{dc}$. This explains the
opposite ring in \eqref{eq:point-end}. At the zero point all variable
actions are zero, so every $c\in D$ is allowed and the endomorphism
ring is $D^{\mathrm{op}}$.

\smallskip\noindent
\emph{Isomorphisms.}
Any $R_n$-isomorphism $D_{\ba}\to D_{\bb}$ is in particular a
left $D$-linear isomorphism, hence is $r_c$ for some $c\in D^\times$.
For each $i$, intertwining the two actions of $x_i$ means
\[
r_cT_i^{\ba}=T_i^{\bb}r_c
\quad\Longleftrightarrow\quad
 a_ic=\sigma(c)b_i.
\]
Solving for $b_i$ and writing $u=c^{-1}$ gives
$b_i=\sigma(u)a_i u^{-1}$. Conversely, if these equations hold, then
$r_{u^{-1}}$ is an invertible left $D$-linear map intertwining every
variable action. It is therefore an $R_n$-isomorphism from $D_\ba$ to $D_\bb$.
\end{proof}

\subsection{One variable and inner twisting}
\begin{lemma}\label{lem:onevar}
Let $A=D[t;\sigma]$. Then:
\begin{enumerate}
	\item Every nonzero proper left ideal $I$ of $A$ can be written uniquely as $I=Ap$ for a monic polynomial $p\in A$. Every simple left $A$-module is isomorphic to $A/Aq$ for some monic irreducible polynomial $q\in A$.
	\item Every simple left $A$-module has finite $D$-dimension.
	\item  All simple $A$-modules have dimension one if and only if $(D,\sigma)$
	is right $\sigma$-algebraically closed.
\end{enumerate}
\end{lemma}
\begin{proof}

Let $p$ be monic of degree $m$, and suppose $f$ has degree $k\ge m$
and leading coefficient $d\in D$. The polynomial $d t^{k-m}p$ has the
same leading term as $f$. Subtracting it strictly decreases the degree.
Repeating this cancellation yields right division by $p$:
\[
f=qp+r,\qquad r=0\text{ or }\deg r<m.
\]
Here the divisor $p$ is on the right, and the remainder $r$ is unique.

Now let $I$ be a nonzero proper left ideal. Choose an element of least degree
and multiply it on the left by the inverse of its leading coefficient
to obtain a monic element $p\in I$. For every $f\in I$, dividing $f$ by $p$ gives
a remainder $r\in I$. The degree minimality  forces $r=0$, so $I=Ap$.
If two monic polynomials $p$ and $q$ generate $I$, each is a left multiple of the
other. Hence they have the same degree. Therefore, the monicity gives $p=q$.

Every nonzero proper left ideal  has the form $Ap$ for
a monic polynomial $p$ of positive degree. $Ap$ is maximal
if and only if $p$ is irreducible. Thus every simple left $A$-module is isomorphic to $A/Ap$ for a monic irreducible polynomial $p$.  
Right division also shows that the residue classes
\[
1+Ap,\ t+Ap,\ \ldots,\ t^{m-1}+Ap
\]
form a left $D$-basis of $A/Ap$, where $m=\deg p$. All irreducible monic polynomial in $A$ has degree 1 if and only if $(D,\sigma)$ is right $\sigma$-algebraically closed. This proves the last equivalence.
\end{proof}

\begin{lemma}[Inner twisting]\label{lem:inner}
Let $u\in D^\times$ and set $\theta=\Int(u^{-1})\circ\sigma$.
The assignments $y_i\mapsto u^{-1}x_i$ define a ring isomorphism
\[
R_n(D,\theta)=D[y_1,\dots,y_n;\theta]\longrightarrow R_n(D,\sigma)=D[x_1,\dots,x_n;\sigma]
\]
that fixes $D$ pointwise. The condition $\Hcond_n$ is invariant under
this replacement of $\sigma$. The same is true of both the all-simple
and the finite-dimensional-simple module conditions in
Theorem~\ref{thm:main}.
\end{lemma}
\begin{proof}
Put $y_i'=u^{-1}x_i$ in $R_n(D,\sigma)$. For $d\in D$,
\[
y_i'd=u^{-1}\sigma(d)x_i
     =\bigl(u^{-1}\sigma(d)u\bigr)y_i'
     =\theta(d)y_i'.
\]
Moreover,
\[
y_i'y_j'=u^{-1}\sigma(u^{-1})x_ix_j
         =u^{-1}\sigma(u^{-1})x_jx_i=y_j'y_i'.
\]
Thus the stated assignments $y_i\mapsto u^{-1}x_i$ define a ring homomorphism
fixing $D$.

Conversely, in $R_n(D,\theta)$ the elements $uy_i$ induce $\sigma$
on $D$, since $u\theta(d)=\sigma(d)u$. They commute with one
another, because
\[
(uy_i)(uy_j)=u\theta(u)y_iy_j
           =u\theta(u)y_jy_i.
\]
They therefore define a homomorphism in the other direction. The two
maps are inverse on $D$ and on each variable, proving the isomorphism.

Transporting a module along an isomorphism preserves simplicity. Since
this isomorphism is the identity on $D$, it also preserves its underlying
left $D$-space and its $D$-dimension.  Applying these observations for $n=1$ and using
Lemma~\ref{lem:onevar} proves invariance of algebraic closure: $R_n(D,\theta)$ is right $\theta$-algebraically closed if and only if $R_n(D,\sigma)$ is right $\sigma$-algebraically closed. 

It remains to compare the division rings in the centralizer condition.
For $v\in D^\times$ and $c\in D$,
\begin{align*}
c\in E_v(D,\theta)
&\Longleftrightarrow u^{-1}\sigma(c)uv=vc\\
&\Longleftrightarrow \sigma(c)(uv)=(uv)c\\
&\Longleftrightarrow c\in E_{uv}(D,\sigma).
\end{align*}
These are equal subrings of $D$. As $v$
ranges over $D^\times$, so does $uv$. Their commuting tuples and
centralizers therefore give exactly the same family of conditions in
$\Hcond_n$.
\end{proof}

\subsection{Scalar extension and variable actions}
We record two elementary observations that will be used repeatedly.

\begin{lemma}[Subspaces under scalar extension]\label{lem:scalar}
Let $C\subseteq D$ be division rings and let $V$ be a left $C$-space.
Identify $V$ with its image in $D\otimes_C V$ under $v\mapsto1\otimes v$.
For every $C$-subspace $W\subseteq V$,
\begin{equation}\label{eq:intersection}
V\cap DW=W.
\end{equation}
In particular, if  $V=\bigoplus_\lambda Cv_\lambda$  and $M=\bigoplus_\lambda Dv_\lambda$, the multiplication map
$D\otimes_C V\to M$ is an isomorphism, and \eqref{eq:intersection}
holds inside $M$.
\end{lemma}
\begin{proof}
Choose a left $C$-basis $\{w_i\}_{i\in I}$ of $W$ and extend it to
a left $C$-basis
\[
\{w_i\}_{i\in I}\ \cup\ \{v_j\}_{j\in J}
\]
of $V$. Tensoring the corresponding direct-sum decomposition with the
right $C$-space $D$ gives
\[
D\otimes_C V
 =\bigoplus_{i\in I}D(1\otimes w_i)
  \ \oplus\ \bigoplus_{j\in J}D(1\otimes v_j).
\]
In particular, these tensors are a left $D$-basis. One may also check
their independence by applying the coefficient maps associated with
the chosen $C$-basis.

An element of $V$ has a unique finite expansion in that $C$-basis.
Its image under $v\mapsto1\otimes v$ has the same coefficients,
now viewed in $D$. The displayed $D$-basis shows that a nonzero
expansion cannot map to zero. This justifies the identification of
$V$ with a subspace of the tensor product.

The $D$-span $DW$ consists exactly of the vectors whose coordinates
at every $1\otimes v_j$ are zero. If a vector originally belonging
to $V$ lies in $DW$, its $C$-coordinates at every $v_j$ must
therefore vanish. It belongs to $W$. The opposite inclusion
$W\subseteq V\cap DW$ is immediate, proving \eqref{eq:intersection}.

In the final assertion, the multiplication map is well-defined because
$(dc)v=d(cv)$ for $d\in D$, $c\in C$, and $v\in V$. It sends the
$D$-basis $\{1\otimes v_\lambda\}$ onto the given $D$-basis
$\{v_\lambda\}$ of $M$, and is therefore an isomorphism. The
intersection identity transfers along this map. All expansions involved
are finite even when the indexing sets are infinite, so no
finite-dimensionality assumption is used.
\end{proof}

\begin{lemma}\label{lem:actions}
On a simple $R_n$-module $M$, each variable action $T_i:M\rightarrow M$, $v\mapsto x_iv$ is either zero
or bijective. If $T_1$ is bijective, then
\[
E_i=T_iT_1^{-1}\qquad(2\le i\le n)
\]
are commuting $R_n$-module endomorphisms of $M$.
\end{lemma}
\begin{proof}
Write $T_i(v)=x_i v$. The defining relations say that each $T_i$ is
additive and $\sigma$-semilinear:
\[
T_i(dv)=\sigma(d)T_i(v)\qquad(d\in D,
\ v\in M).
\]
Its kernel and its image 
 are stable under every $T_j$: commutativity gives
$T_iT_j(v)=T_jT_i(v)$ for the kernel and
$T_j(T_i(v))=T_i(T_j(v))$ for the image. They are consequently
$R_n$-submodules of $M$. If $T_i\ne0$, its kernel is not all of $M$ and
its image is nonzero. Simplicity forces the kernel to be zero and the
image to be $M$. This proves that $T_i$ is bijective.

Suppose now that $T_1$ is bijective. Its inverse is
$\sigma^{-1}$-semilinear, as can be checked by applying $T_1$ to both sides of
\[
T_1^{-1}(dv)=\sigma^{-1}(d)T_1^{-1}(v).
\]
Therefore
\[
E_i(dv)=T_i\bigl(\sigma^{-1}(d)T_1^{-1}(v)\bigr)
       =dE_i(v),
\]
so $E_i$ is $D$-linear. The equality $T_1T_j=T_jT_1$ implies
$T_1^{-1}T_j=T_jT_1^{-1}$. Hence each $E_i=T_iT_1^{-1}$ commutes
with every $T_j$. Together with $D$-linearity, this is exactly what
is needed for $E_i$ to be an $R_n$-module endomorphism. Finally,
\[
E_iE_j=T_iT_jT_1^{-2}=T_jT_iT_1^{-2}=E_jE_i.
\]
%Notice that the individual maps $T_i$ need not be $D$-linear; the  cancellation of the two semilinearities is what makes their ratios module endomorphisms.
\end{proof}

\section{Contraction in an auxiliary central variable}\label{sec:contraction}

 The argument below follows Amitsur--Small type: a leading-term reduction
turns a suitable central nonunit into a nonzero contraction
\cite{AmitsurSmall,ArySkew}. Here the required central polynomials are
chosen over the prime field. This removes any need to assume that
$Z(D)^\sigma$ is infinite.

\begin{lemma}[Simultaneous avoidance]\label{lem:avoid}
Let $P$ be the prime field of $D$, let $B=D[z]$ with $z$ central,
and let $h_1,\ldots,h_s\in B$ be nonzero. There exists a nonconstant
$r\in P[z]$ such that, for every $B/Bh_j$, the left multiplication action $f+Bh_j\mapsto rf+Bh_j$ is bijective on $B/Bh_j$.

\end{lemma}
\begin{proof}
For each $j$, right division in the central polynomial ring $B=D[z]$
shows that $V_j=B/Bh_j$ has left $D$-dimension $\deg h_j$.
Every polynomial $r\in P[z]$ is central in $B$, so left
multiplication by $r$ induces a well-defined $D$-linear map on $V_j$.

Fix one of these spaces, say $V=V_j$. We claim that only finitely many
monic irreducibles $r\in P[z]$ have a nonzero kernel on $V$.
To see this, take distinct monic irreducibles $r_1,\ldots,r_k$ and
vectors $v_i\in\ker(r_i|_V)$ such that $\sum_{i=1}^k v_i=0$. For a fixed
$i$, put $q_i=\prod_{\ell\ne i}r_\ell$. Then $q_iv_l=0$ for all $l\ne i$, so
$q_iv_i=q_i(\sum_{j=1}^k v_j)=0$.
The polynomials $r_i$ and $q_i$ are relatively prime in the ordinary commutative polynomial  ring $P[z]$. Thus there are $s_i,t_i\in P[z]$ such that
\[
s_ir_i+t_iq_i=1.
\]
Acting on $v_i$ gives $v_i=0$. This proves that the subspaces
$\ker(r_i|_V)$ form a direct sum.

Every nonzero kernel contributes at least one to the $D$-dimension of
that direct sum. Consequently at most $\dim_DV$ monic irreducibles
can have nonzero kernel on $V$.  Applying this argument to the family $\{V_j\}_{j=1}^s$ show that there are only finitely many monic irreducible polynomials $r_i$ that has at least a nonzero kernel. The fact that  there are infinitely many monic irreducible polynomials over any field completes our proof.
\end{proof}

\begin{theorem}[Central-variable contraction]\label{thm:contraction}
Let $z$ be central over $R_n(D,\sigma)$. Every maximal left ideal
$L$ of $R_n[z]$ satisfies
\begin{equation}\label{eq:contraction}
L\cap D[z]\ne0.
\end{equation}
\end{theorem}
\begin{proof}
Set $B=D[z]$ and extend $\sigma$ to $B$ by fixing $z$. Then
$R_n[z]=B[x_1,\ldots,x_n;\sigma]$. The polynomial ring $P[z]$, where $P$ is the prime field of $D$, belongs to the center of $R_n[z]$. The aim is to choose a central
polynomial $r\in B$ for which leading coefficients in $L$ are
saturated under multiplication by $r$. Maximality will then reduce
an equation modulo $L$ to one involving only $z$.

\smallskip\noindent
\emph{Step 1: The leading coefficient ideals.} 
Choose a monomial order in the $x$-variables: it is a well-order,
has $1$ as its least monomial, and is preserved by multiplication
by a monomial. Write $\operatorname{lc}(f)$ for leading coefficient of the leading term of $f\in R_n[z]$: we always take leading terms with respect to the $x$-variables,
 regarding elements of $B$ as coefficients. 
For $\alpha\in\N^n$, define $$I_\alpha=\{0\}\cup\{\operatorname{lc}(f): f\in L\text{ whose leading monomial is }x^\alpha \}.$$ 
The set $I_\alpha$ is a left ideal of $B$. The left
multiplication by $x^\beta$ gives
\[
\sigma^{|\beta|}(I_\alpha)\subseteq I_{\alpha+\beta}.
\]
The twists can be removed by setting
$J_\alpha=\sigma^{-|\alpha|}(I_\alpha)$. Applying the appropriate
inverse power of $\sigma$ to this inclusion yields
\[
J_\alpha\subseteq J_{\alpha+\beta}
\qquad(\alpha,\beta\in\N^n).
\]
Each $J_\alpha$ is again a left ideal of $B$.

\smallskip\noindent
\emph{Step 2: Choose one central polynomial for all leading positions.}
The set
\[
\mathcal S=\{\alpha\in\N^n:J_\alpha\ne0\}
\]
is upward closed for the coordinatewise order. Dickson's lemma says
that such a set has finitely many minimal elements and that every
member dominates one of them. Write these minimal elements as
$\beta_1,\ldots,\beta_s$, and choose
$0\ne h_j\in J_{\beta_j}$ for each $j=1,\dots,s$. By Lemma~\ref{lem:avoid}, there is a
nonconstant $r\in P[z]$ acting bijectively on every $B/Bh_j$.
%If $\mathcal S$ is empty, choose any nonconstant $r\in P[z]$.

%Every automorphism of $D$ fixes its prime field pointwise. Thus $\sigma(r)=r$. Since $P[z]\subseteq Z(B)$, the relations with the $x_i$ show that $r$ is central in the whole ring $R_n[z]$.

\smallskip\noindent
\emph{Step 3: Prove saturation of the leading-coefficient ideals.}
First let $\alpha\in\mathcal S$, and choose $j$ such that
$\beta_j\le\alpha$ coordinatewise. Then
\[
Bh_j\subseteq J_{\beta_j}\subseteq J_\alpha\subseteq B.
\]
Consider the finite-dimensional $D$-space $Q=B/Bh_j$ and its
subspace $U=J_\alpha/Bh_j$. Multiplication by $r$ is bijective on
$Q$ and preserves $U$. Its restriction to $U$ is injective and hence
surjective, so $rU=U$. It follows that the induced map on $Q/U$ is
injective: if $r\xi\in U$, choose $\eta\in U$ with
$r\eta=r\xi$ and use injectivity on $Q$ to obtain $\xi=\eta$.
Since $Q/U\simeq B/J_\alpha$, we have proved
\[
br\in J_\alpha\quad\Longrightarrow\quad b\in J_\alpha.
\]
It follows that
\begin{equation}\label{eq:saturation}
br\in I_\alpha\quad\Longrightarrow\quad b\in I_\alpha.
\end{equation}
If $\alpha\notin\mathcal S$, then $I_\alpha=0$, and the same
implication follows from $r\ne0$. Thus \eqref{eq:saturation} holds for every $\alpha$.

\smallskip\noindent
\emph{Step 4: Eliminate the skew variables.}
If $r\in L$, then $r$ itself is a nonzero element of $L\cap B$.
Otherwise, $L+R_n[z]r=R_n[z]$, so there is a
$0\ne q\in R_n[z]$ with $qr-1\in L$.  Among all
such $q$, choose one whose leading $x$-monomial is least.

Suppose that $\lm(q)=x^\alpha\ne1$, and write $b=\lc(q)\in D[z]$.
Because $r$ is in $P[z]\subset D[z]$,
$qr-1$ still has leading monomial $x^\alpha$, with coefficient
$br$. Therefore $br\in I_\alpha$, and saturation gives
$b\in I_\alpha$. By definition of $I_\alpha$, choose $h\in L$
with leading term $bx^\alpha$. Centrality is crucial here: although
$L$ is only a left ideal, it implies $hr=rh\in L$. Consequently,
\[
(q-h)r-1=(qr-1)-hr\in L.
\]
If $q-h=0$, this puts $-1$ in $L$, which is impossible. Otherwise
the leading terms of $q$ and $h$ cancel, so $q-h$ has a strictly
smaller leading monomial. This contradicts the choice of $q$.

We conclude that $\lm(q)=1$, hence $q\in B$. Now $0\ne qr-1\in L\cap B$.
\end{proof}

\begin{corollary}[Local torsion for endomorphisms]\label{cor:local}
Let $M\ne 0$ be a simple $R_n$-module and let $E\in\End_{R_n}(M)$.
For every $v\in M$, the space
\[
D[E]v=\Span_D\{v,Ev,E^2v,\ldots\}
\]
is finite-dimensional over $D$. { Regarding $M$ as a left
$D[z]$-module, with $z$ acting as $E$, then $M$ contains a finite-dimensional simple
$D[z]$-submodule.}
\end{corollary}
\begin{proof}
Because $E$ is an $R_n$-module endomorphism, it commutes with the
action of every element of $R_n$. We may therefore extend the action
to $R_n[z]$ by letting the central variable $z$ act as $E$.
An $R_n[z]$-submodule is in particular an $R_n$-submodule, so $M$
remains simple through this extended action.

Fix $v\ne0$. Since $M$ is simple, the map $R_n[z]\to M$ sending $f$ to $fv$ is
surjective and its kernel $\Ann_{R_n[z]}(v)$ is a maximal left ideal. By Theorem~\ref{thm:contraction}, this kernel
contains a nonzero polynomial
\[
q(z)=q_0+q_1z+\cdots+q_mz^m\in D[z],\qquad m\ge 1,~~q_m\ne0.
\]
The equation $q(E)v=0$ can be re-written as
\[
E^mv=-\sum_{j=0}^{m-1}q_m^{-1}q_jE^jv.
\]
It follows that $D[E]v\subseteq \sum_{i=0}^{m-1}DE^iv$, which is finite-dimensional over $D$.

The space $D[E]v$ is also a nonzero $D[z]$-submodule,
since it is stable under $D$ and under $E$. Among its nonzero
$D[z]$-submodules, choose one of least positive $D$-dimension. A
proper nonzero submodule would have smaller dimension, so the chosen
one is simple.
\end{proof}

\section{The finite-dimensional criterion}\label{sec:finite}

The following lemma follows from
\cite[Theorem~4.4, implication $(3)\Rightarrow(2)$]{AryCentral}.

\begin{lemma}\label{lem:centralfinite}
 Let $E$ be a division ring and let $m\ge1$. Suppose
$\Cen_E(b_1,\ldots,b_{m-1})$ is right algebraically closed for every
pairwise commuting $(m-1)$-tuple $(b_1,\ldots,b_{m-1})\in E^{m-1}$. Then every
finite-dimensional simple left $E[z_1,\ldots,z_m]$-module has
$E$-dimension one.
\end{lemma}

\begin{proposition}\label{prop:finiteskew}
If $\Hcond_n(D,\sigma)$ holds, every finite-dimensional simple
$R_n$-module has $D$-dimension one.
\end{proposition}
\begin{proof}
The case $n=1$ is Lemma~\ref{lem:onevar}. Let $n\ge2$, and let
$M\ne0$ be a finite-dimensional simple $R_n$-module with variable actions
$T_1,\ldots,T_n$. If all variable actions are zero, every left
$D$-subspace of $M$ is an $R_n$-submodule; simplicity then gives
$\dim_DM=1$. Otherwise, after relabeling, Lemma~\ref{lem:actions}
allows us to take $T_1$ bijective. Set
$E_i=T_iT_1^{-1}$ for $2\le i\le n$.

\smallskip\noindent
\emph{Step 1: Restrict to one skew variable.}
Let { $A=D[x_1;\sigma]$} act on $M$ through $T_1$. The first part of
$\Hcond_n$ and Lemma~\ref{lem:onevar} imply that every simple
$A$-module has $D$-dimension one.
Let $N$ be a nonzero $A$-submodule of $M$ of minimal dimension over $D$. For any $0\ne v\in N$, then $Av\subseteq M$ is an $A$-submodule of $M$. So $\dim_D(Av)=\dim_D(N)$. It follows that $Av=N$. This shows that $N$ is a simple $A$-submodule of $M$.

%Since $M$ is finite-dimensional, its restriction to $A$ contains a simple $A$-submodule, obtained by minimizing the positive $D$-dimension of a submodule.

Each $E_i$ is $A$-linear by Lemma~\ref{lem:actions}. It maps a
simple $A$-submodule of $M$ to zero or to an isomorphic simple submodule.
The $A$-socle of $M$ is therefore stable under every $E_i$; it is also
stable under $T_1$, as any $A$-submodule is. Thus it is stable under
$T_i=E_iT_1$ for all $i$, and is a nonzero $R_n$-submodule. It
must equal $M$. Grouping its simple summands into isotypic components
$$M=\oplus_\lambda M_\lambda^{s_\lambda}.$$

Each isotypic component $M_\lambda^{s_\lambda}$ is stable under $E_i,T_1,$ and $T_i$ for all $i$. Hence it is a $R_n$-submodule of $M$. 
 Simplicity of $M$ thus
forces a single $A$-submodule type:
$$M\simeq M_0^s\qquad\text{as $A$-module}$$
for one simple $A$-module $M_0=D$.

Choose a nonzero generator \(v_1\) of one copy and write $T_1v_1=av_1$
for some \( a\in D\). Because \(T_1\) is injective, $a\ne0.$ 
For each other copy of \(M_0\), choose an \(A\)-isomorphism from the first copy and let \(v_j\) be the image of \(v_1\). The isomorphism is \(D\)-linear and commutes with \(T_1\), so
\[ T_1v_j=av_j. \]
Thus we obtain a $D$-basis $v_1,\ldots,v_s$ of $M$ such that
\[
T_1v_j=av_j\qquad(1\le j\le s).
\]

\smallskip\noindent
\emph{Step 2: Pass to an eigenspace over $E_a$.}
Define $V=\{v\in M:T_1v=av\}$. For $v=\sum_jd_jv_j$,
$D$-semilinearity of $T_1$ gives
\[
T_1v=\sum_j\sigma(d_j)av_j,
\qquad av=\sum_jad_jv_j.
\]
Therefore $v\in V$ if and only if $\sigma(d_j)a=ad_j$ for every
$j$, or equivalently $d_j\in E_a=D^{\theta_a}$. Consequently,
\[
V=\bigoplus_{j=1}^s E_av_j,
\qquad M\simeq D\otimes_{E_a}V.
\]
For $v\in V$, the $D$-linearity and commutation properties of $E_i$
give
\[
T_1(E_iv)=E_i(T_1v)=E_i(av)=aE_iv.
\]
Thus each $E_i$ preserves $V$. Their restrictions are commuting
$E_a$-linear operators, so they define a module structure over the
central polynomial ring $E_a[z_2,\ldots,z_n]$, with $z_i$ acting as
$E_i$.

\smallskip\noindent
\emph{Step 3: Apply the finite-dimensional criterion.}
We first check that $V$ is a simple $E_a[z_2,\ldots,z_n]$-module. For a submodule
$0\ne W\subseteq V$, take its $D$-span in $M$:
\[
DW=
\left\{
\sum_{k=1}^r d_kw_k:
d_k\in D,\ w_k\in W
\right\}\subseteq M.
\]
Then for $d\in D$ and $w\in W$,
\[
T_1(dw)=\sigma(d)aw\in DW,
\qquad
T_i(dw)=\sigma(d)E_i(aw)=\sigma(d)aE_iw\in DW.
\]
Hence $DW$ is a nonzero $R_n$-submodule of $M$, so $DW=M$.
Lemma~\ref{lem:scalar} gives $$W=V\cap DW=V\cap M=V.$$
So $V$ is a simple $E_a[z_2,\ldots,z_n]$-module.  By the second part of $\Hcond_n$, we can apply
Lemma~\ref{lem:centralfinite} for the division ring $E_a$ to obtain that $\dim_{E_a}V=1$. It follows that $$\dim_DM=\dim_D( D\otimes_{E_a}E_a)=1.$$
\end{proof}

\section{Finite orbits and descent to a centralizer}\label{sec:descent}

\subsection{Splitting a commuting endomorphism}
The next lemma is the step that allows the argument to leave the
finite-dimensional category. Its hypothesis concerns only
finite-dimensional simple modules in two skew variables.

\begin{lemma}[Finite-orbit]\label{lem:splitting}
Suppose that every finite-dimensional simple $R_2(D,\sigma)$-module has
$D$-dimension one. Let $M$ be a simple $R_n(D,\sigma)$-module with
at least one nonzero variable action, and let $E\in\End_{R_n}(M)$. Regard $M$
as a module over $B=D[z]$ with $z$ acting as $E$, then $M$ is
a direct sum of copies of a one-dimensional simple $B$-module $L$.
There are $b\in D$ and  $u\in D^\times$ satisfying 
\begin{equation}\label{eq:fixedb}
\sigma(b)u=ub,
\end{equation}
such that $E$ act on $L$
as the right multiplication by $b$.
\end{lemma}
\begin{proof}
Relabel the variables so that $T_1\ne0$. By Lemma~\ref{lem:actions},
$T_1$ is bijective. Throughout the proof, $B=D[z]$ acts through
$z=E$. The original module $M$ is not assumed finite-dimensional;
we will construct a finite-dimensional auxiliary module to which the
hypothesis can be applied.

\smallskip\noindent
\emph{Step 1: Show that $M$ is semisimple over $B$.}
Corollary~\ref{cor:local} supplies a finite-dimensional simple
$B$-submodule $W\subseteq M$. For $\alpha\in\N^n$, write
$T^\alpha=T_1^{\alpha_1}\cdots T_n^{\alpha_n}$ and put
$d=|\alpha|$. If $T^\alpha W\ne0$, none of the variable actions
occurring in this product is zero. Each is bijective by
Lemma~\ref{lem:actions}, so $T^\alpha$ is a bijective
$\sigma^d$-semilinear map on $M$. It commutes with $E$.
For $f\in B$, extending $\sigma$ by $\sigma(z)=z$, we have
\[
T^\alpha(fw)=\sigma^d(f)T^\alpha(w),\quad \forall w\in W.
\]
This identity shows that $T^\alpha W$ is a $B$-submodule and that
any of its $B$-submodules pulls back to a $B$-submodule of $W$.
Thus $T^\alpha W$ is simple.

The $R_n$-submodule generated by $W$ is nonzero and hence equals
$M$. Since each $T^\alpha W$ is already a left $D$-space,
\begin{equation}\label{eq:sum-simple}
M=R_nW=\sum_{\alpha\in\N^n}T^\alpha W.
\end{equation}
We have expressed $M$ as a sum of simple $B$-submodules. Zorn's lemma gives a maximal collection of these
submodules with direct sum, this maximal direct sum exhausts $M$: if there exists any remaining simple
submodule not contained in this direct sum, then it would have zero intersection
with the sum and could be added to the sum, this contradicts the maximality of the collection. We can rewrite (\ref{eq:sum-simple}) as
\begin{equation}\label{eq:directsum-simplemodules}
M=\bigoplus_{\alpha\in\Gamma}T^\alpha W, \qquad T^\alpha W\ne0,
\end{equation}
where $\Gamma\subseteq \mathbb{N}^n$.

\smallskip\noindent
\emph{Step 2: Prove that the simple types form a finite orbit.} For each $\alpha\in\Gamma$, set $d=|\alpha|$, the surjective map \[
h_\alpha=T^\alpha T_1^{-d}:
T_1^dW\longrightarrow T^\alpha W.
\] is a $B$-module homomorphism. The source \(T_1^dW\) is simple over \(B\), so it is injective as well. We obtain that
\begin{equation}\label{eq:T-alpha-iso}
    T^\alpha W\cong_B T_1^dW,\qquad d=|\alpha|.
\end{equation}

Now let \(\mathcal S\) be the set of isomorphism classes of simple \(B\)-submodules of \(M\), and put
\[
\tau=[W]\in\mathcal S.
\]
Here \([W]\) denotes the isomorphism type of \(W\), so
\[
[U]=[W]\quad\Longleftrightarrow\quad U\cong_BW.
\]
Define the operation
$$F:\mathcal S\longrightarrow \mathcal S,\qquad [U]\longmapsto[T_1U].$$
This is well-defined: if
$f:U\to U'$ is a $B$-isomorphism, then $T_1fT_1^{-1}$ is a
$B$-isomorphism $T_1U\to T_1U'$. Moreover, $F$ is bijective because the inverse  is induced
by $T_1^{-1}$. In particular,
\[
F^d(\tau)=[T_1^dW],
\qquad
F^{-1}(\tau)=[T_1^{-1}W].
\]

We now show that the orbit of \(\tau\) under \(F\) is finite.

 Since  $T_1^{-1}W$ is a simple $B$-submodule of $M$, the direct sum (\ref{eq:directsum-simplemodules}) gives $T_1^{-1}W\cong T^\alpha W$ for some $\alpha\in\Gamma$. It follows, by equation (\ref{eq:T-alpha-iso}), that $T_1^{-1}W\cong T_1^dW$ where $d=|\alpha|$. Hence
\[
F^{-1}(\tau)=F^d(\tau).
\]
Applying $F$ gives $F^{d+1}(\tau)=\tau$. The orbit is consequently
finite. Let $\ell\ge1$ be its length, and choose a $B$-isomorphism
$\phi:T_1^\ell W\to W$.

\smallskip\noindent
\emph{Step 3: Construct an $R_2$-module.}
Form the \emph{external} direct sum
\begin{equation}\label{eq:directsum-simplemodules2}
    N=W\oplus T_1W\oplus\cdots\oplus T_1^{\ell-1}W.
\end{equation}
The use of an external sum does not require these submodules to be
independent inside $M$. Each summand has the same finite
$D$-dimension as $W$, so $\dim_DN=\ell\dim_DW$.
Define $T$ from a summand to the next by $T_1$, and from the last
summand back to the first by $\phi T_1$. It is a bijective
$\sigma$-semilinear map. Let $E_N$ act on each summand by the original
operator $E$. The operators $T$ and $E_N$ commute. Therefore $T$ and $E_NT$ are
commuting $\sigma$-semilinear operators, making $N$ an $R_2(D,\sigma)$-module.

\smallskip\noindent
\emph{Step 4: Use the finite-dimensional hypothesis.}
Choose a nonzero $R_2$-submodule $L\subseteq N$ of least
$D$-dimension. It is simple, and the hypothesis gives $\dim_DL=1$.
Since $T$ is injective, it restricts to a bijection
of $L$: $T(L)=L$. In particular, $T^{-1}(L)=L$. Hence $E_N=(E_NT)T^{-1}$ also preserves $L$.
Thus $L$ is also a one-dimensional simple $B$-submodule of $N$.
Since every summand in the direct sum (\ref{eq:directsum-simplemodules2}) is a simple $B$-module, $L$ is isomorphic to $T_1^jW$ for some $j$. Let $h:L\rightarrow T_1^jW$ be a $B$-module isomorphism. Then \[
g=(T|_L)\circ h^{-1}\circ
\bigl(T_1^{-1}|_{T_1^{j+1}W}\bigr):T_1^{j+1}W\longrightarrow L.
\]
is also a $B$-module isomorphism. Therefore,
\[
T_1^{j+1}W\cong T_1^jW,\quad\text{or }~~~F^{j+1}(\tau)=F^j(\tau).
\]
Applying $F^{-j}$ gives $F(\tau)=\tau.$ Consequently, $\ell=1$ and \(\dim_DW=1\); hence \(N=W=L\) after identification. The equations (\ref{eq:T-alpha-iso})  now show that all simple $B$-summands of $M$ in (\ref{eq:directsum-simplemodules2}) are therefore one-dimensional and
have the same type: it is the one dimensional simple $B$-module $W$.

Finally, choose $0\ne v\in L$ and write
\[
Tv=uv,\qquad E_Nv=bv
\]
with $u\in D^\times$ and $b\in D$. Evaluating $TE_N=E_NT$ at
$v$ gives
\[
\sigma(b)uv=ubv,
\]
so $\sigma(b)u=ub$. $L$ is generated by $v$, the action of \(z\) is
\[
\begin{aligned}
z\cdot(dv)
&=E_N(dv)=dE_Nv=dbv.
\end{aligned}
\]
Therefore the corresponding operator on \(D\) is $d\longmapsto db.$
That is, \(z\) acts by right multiplication by \(b\).
\end{proof}

\subsection{Heredity of the closure condition}
\begin{lemma}[Centralizer descent]\label{lem:heredity}
Assume $\Hcond_n(D,\sigma)$ with $n\ge2$. Suppose $b\in D$ and
$u\in D^\times$ satisfy \eqref{eq:fixedb}, and put
\[
C=\Cen_D(b),
\qquad
\theta=\Int(u^{-1})\circ\sigma.
\]
Then $\theta(b)=b$, the automorphism $\theta$ restricts to $C$, and
$\Hcond_{n-1}(C,\theta|_C)$ holds.
\end{lemma}
\begin{proof}
The relation $\sigma(b)u=ub$ gives
$\theta(b)=u^{-1}\sigma(b)u=b$. If $c\in C$, applying $\theta$ to
$cb=bc$ gives $\theta(c)b=b\theta(c)$. Thus $\theta(C)\subseteq C$.
The same argument applies for $\theta^{-1}$ , so the restriction
$\theta|_C$ is an automorphism. By Lemma~\ref{lem:inner}, the
hypothesis also gives $\Hcond_n(D,\theta)$.

\smallskip\noindent
\emph{Step 1: Prove right $\theta$-algebraic closure of $C$.}
Suppose this closure fails. By Lemma~\ref{lem:onevar}, there is a
monic irreducible polynomial
\[
p=t^m+p_{m-1}t^{m-1}+\cdots+p_0\in C[t;\theta],
\qquad m\ge2.
\]
Its constant coefficient $p_0$ is nonzero. Otherwise $p$ would have
the proper right factor $t$, contradicting irreducibility.
Set $A=D[t;\theta]$ and consider the left  $A$-module
\[
M=A/Ap.
\]
This is an $m$-dimensional left $D$-space. Every coset has a unique
representative of degree less than $m$, by right division by the monic
polynomial $p$.

Let $T$ be the action of $t$ on $M$. It is
$\theta$-semilinear and bijective. Next define $B(\bar q)=\overline{qb}$ for $\bar{q}\in A/Ap$. We check the well-definedness of this right multiplication by $\bar{b}$: since $\theta(b)=b$,
we have $tb=bt$; and every coefficient of $p$ commutes with $b$; thus $pb=bp$, and for any $h\in A$,
\[
(hp)b=h(pb)=hb\,p\in Ap.
\]
Therefore $B$ is well-defined. The operators $T$ and $BT$ are commuting
$\theta$-semilinear maps, so they give $M$ the structure of an
$R_2(D,\theta)$-module.

Now choose a nonzero simple \(R_2(D,\theta)\)-submodule $W\subseteq M$. 
Such a submodule exists by minimizing positive \(D\)-dimension. Since \(\mathcal H_n(D,\theta)\) implies \(\mathcal H_2(D,\theta)\), Proposition 4.2 gives
\[
W=D\bar q
\]
for a monic reduced representative \(q\in D[t;\theta]\) with $k=\deg q<m$.

Because \(T\) is injective, \(T(W)=W\). Hence
$B=(BT)T^{-1}$
also preserves \(W\). Consequently, for some \(c,d\in D\),
\[
T(\bar q)=c\bar q,\qquad B(\bar q)=d\bar q.
\]
or
\[
tq-cq\in Ap,\qquad qb-dq\in Ap.
\]
The second polynomial has degree less than \(m\), so it is zero. Since \(\theta(b)=b\), the variable \(t\) commutes with \(b\). Comparing leading coefficients in \(qb=dq\) therefore gives \(d=b\), and then $qb=bq$.
Thus every coefficient of \(q\) belongs to \(C=\mathrm C_D(b)\): $q\in C[t;\theta]$.
The polynomial \((t-c)q\) is monic, nonzero, has degree at most \(m\), and belongs to \(Ap\). Therefore
\[
k=m-1,\qquad p=(t-c)q.
\]
Writing
\[
q=t^{m-1}+q_{m-2}t^{m-2}+\cdots,
\]
comparison of the coefficient of \(t^{m-1}\) gives
\[
{c=\theta(q_{m-2})-p_{m-1}\in C.}
\]
Both factors belong to \(C[t;\theta]\) and have positive degree, contradicting irreducibility of \(p\). This proves the required closure of \(C\).

\smallskip\noindent
\emph{Step 2:Transfer the centralizer conditions.}
For $n=2$, the conclusion is just the one-variable closure proved
above. Let $n\ge3$, and fix $v\in C^\times$. On $D$, put
\[
\psi=\Int(v^{-1})\circ\theta
     =\Int((uv)^{-1})\circ\sigma.
\]
Because $v$ commutes with $b$ and $\theta(b)=b$, we also have
$\psi(b)=b$. The fixed elements of $\psi$ in $C$ are therefore
\[
E_v(C,\theta)
 =C\cap D^\psi
 =\Cen_{E_{uv}(D,\sigma)}(b).
\]
Take a pairwise commuting $(n-3)$-tuple $\bc$ in this ring. Each
entry of $\bc$ commutes with $b$, and all entries lie in
$E_{uv}(D,\sigma)$. Hence $(b,\bc)$ is a commuting tuple of length
$n-2$ in that division ring. Directly from the definitions,
\[
\Cen_{E_v(C,\theta)}(\bc)
 =\Cen_{E_{uv}(D,\sigma)}(b,\bc).
\]
The right-hand side is right algebraically closed by $\Hcond_n$.
As $v$ and $\bc$ were arbitrary, this is exactly the second
condition of $\Hcond_{n-1}(C,\theta|_C)$.
\end{proof}

\section{The dimension of simple $R_n$-module}\label{sec:proof}

\subsection{Sufficiency for arbitrary simple modules}
\begin{proposition}\label{prop:sufficiency}
If $\Hcond_n(D,\sigma)$ holds, every simple $R_n$-module has
$D$-dimension one.
\end{proposition}
\begin{proof}
We use induction on $n$, uniformly over all division rings and all
automorphisms. This formulation allows the coefficient division ring
to change at the induction step. For $n=1$, the assertion is
Lemma~\ref{lem:onevar}.

Let $n\ge2$, assume the result for $n-1$, and let $M$ be a simple
$R_n(D,\sigma)$-module. No restriction is placed on $\dim_DM$.
If all variables act as zero, any nonzero vector spans a nonzero
$R_n$-submodule over $D$, so simplicity gives $\dim_DM=1$.
Otherwise relabel so that $T_1$ is bijective, and put
\[
E=T_nT_1^{-1}.
\]
By Lemma~\ref{lem:actions}, this is a $D$-linear endomorphism of the
$R_n$-module $M$.

\smallskip\noindent
\emph{Step 1: Prepresentation through a set of eigenvectors.}
The hypothesis $\Hcond_n$ implies $\Hcond_2$. Hence applying 
Proposition~\ref{prop:finiteskew} gives that every finite-demensional simple $R_2$ module has $D$-dimension one. Consequently this verifies the
finite-dimensional hypothesis of Lemma~\ref{lem:splitting}. 
Applying Lemma~\ref{lem:splitting} with  $E\in\operatorname{End}_{R_n}(M)$ shows that, over $D[z]$ with $z=E$,
the module $M$ is a direct sum of copies of a one-dimesional simpe $D[z]$-module $L_0=Dv_0$:
$$M=\bigoplus_{\lambda\in \Lambda}L_\lambda,$$
where there is a $D[z]$-isomorphism $\phi_\lambda: L_0\cong L_\lambda$ for every $\lambda\ne0$. The set $\Lambda$ may be infinite. Moreover, there are elements  $b\in D$ and 
$u\in D^\times$ satisfying $\sigma(b)u=ub$ such that  $$z(dv_0)=Edv_0=dbv_0,\quad d\in D.$$
In particular, $Ev_0=bv_0.$ Set $v_\lambda=\phi_\lambda(v_0)$, we obtain
a left $D$-basis $\{v_\lambda\}_{\lambda\in\Lambda}$ satisfying
\[
Ev_\lambda=bv_\lambda\qquad(\lambda\in\Lambda).
\]

Set
\[
C=\Cen_D(b),\qquad
\theta=\Int(u^{-1})\circ\sigma,\qquad
V=\{v\in M:Ev=bv\}.
\]
The equation $\sigma(b)u=ub$ gives $\theta(b)=b$ and makes $C$
stable under $\theta$. To determine $V$, write a vector $v$ as a finite
sum $v=\sum_\lambda d_\lambda v_\lambda$. The $D$-linearity of
$E$ gives
\[
Ev=\sum_\lambda d_\lambda b v_\lambda,
\qquad
bv=\sum_\lambda b d_\lambda v_\lambda.
\]
Thus $v\in V$ if and only if every $d_\lambda$ commutes with $b$.
We have obtained
\[
M=\bigoplus_{\lambda\in\Lambda}Dv_\lambda,
\qquad
V=\bigoplus_{\lambda\in\Lambda}Cv_\lambda.
\]
Lemma~\ref{lem:scalar} now yields
\begin{equation}\label{eq:descent-tensor}
D\otimes_CV\xrightarrow{\ \sim\ }M,
\qquad d\otimes v\longmapsto dv.
\end{equation}
In particular, this scalar-extension description is valid before we
know that either of these spaces is finite-dimensional.

\smallskip\noindent
\emph{Step 2: Descend the variable actions from $M$ to $V$.}
Define the action $U_i=u^{-1}T_i$ on $M$. For $v\in V$, $D$-linearity of
$E$ and its commutation with $T_i$ give
\begin{align*}
E(U_iv)
 &=u^{-1}T_i(Ev)
  =u^{-1}\sigma(b)T_iv\\
 &=b u^{-1}T_iv=bU_iv.
\end{align*}
Hence every $U_i$ preserves $V$. For $c\in C$ and $v\in V$,
\[
U_i(cv)=u^{-1}\sigma(c)T_iv
       =\bigl(u^{-1}\sigma(c)u\bigr)U_iv
       =\theta(c)U_iv.
\]
Their restrictions are therefore $\theta|_C$-semilinear.
They commute, since the calculation on the ambient space $M$ gives
\[
U_iU_j=u^{-1}\sigma(u^{-1})T_iT_j
      =u^{-1}\sigma(u^{-1})T_jT_i=U_jU_i.
\]
Finally, $T_n=ET_1$ and $D$-linearity of $E$ give
$U_n=EU_1$ on $M$. Because $U_1$ preserves $V$, this becomes
\[
U_n=bU_1\quad\text{on }V.
\]
The element $b$ belongs to $Z(C)$: every element of $C$ commutes
with it. Also $\theta(b)=b$.
These operators define an \(R_{n-1}(C,\theta|_C)\)-module structure on \(V\), with \(x_i\) acting as \(U_i\).

% The first $n-1$ operators $U_1,\dots,U_{n-1}$ give $V$ an $R_{n-1}(C,\theta|_C)\equiv C[U_1,\dots,U_{n-1};\theta]$-module structure. %and the last operator is determined by the first through the displayed identity.

\smallskip\noindent
\emph{Step 3: Prove that the descended module is simple.}
Let $0\ne W\subseteq V$ be a submodule for the first $n-1$
operators. It is a left $C$-subspace and is stable under $U_n$ as
well, since $U_n=bU_1$ and $b\in C$. For every $d\in D$,
$w\in W$, and $1\le i\le n$,
\[
T_i(dw)=\sigma(d)T_iw=\sigma(d)uU_iw\in DW.
\]
Thus $DW$ is a nonzero $R_n$-submodule of $M$, and simplicity gives
$DW=M$. By \eqref{eq:descent-tensor} and the intersection formula
in Lemma~\ref{lem:scalar},
\[
W=V\cap DW=V.
\]
So $V$ is simple over $R_{n-1}(C,\theta|_C)$.

\smallskip\noindent
\emph{Step 4: Apply induction over the smaller division ring.}
Lemma~\ref{lem:heredity} gives
$\Hcond_{n-1}(C,\theta|_C)$. The induction hypothesis therefore
applies to the simple module $V$ and yields $\dim_CV=1$.
Equation~\eqref{eq:descent-tensor} then gives $\dim_DM=1$. % The only finite-dimensional modules used before this conclusion were the auxiliary modules in Lemma~\ref{lem:splitting}; finite-dimensionality of the original module was never assumed.
\end{proof}

\subsection{Necessity from finite-dimensional modules}
\begin{proposition}\label{prop:necessity}
If every finite-dimensional simple $R_n(D,\sigma)$-module has
$D$-dimension one, then $\Hcond_n(D,\sigma)$ holds.
\end{proposition}
\begin{proof}
First consider a simple $D[t;\sigma]$-module $N$.
Lemma~\ref{lem:onevar} makes it finite-dimensional over $D$.
Let $x_1$ act as $t$ and let $x_2,\ldots,x_n$ act as zero.
These actions satisfy all defining relations of $R_n$, and the
$R_n$-submodules are exactly the $D[t;\sigma]$-submodules.
Thus $N$ becomes a finite-dimensional simple $R_n$-module, so the
hypothesis gives $\dim_DN=1$. Lemma~\ref{lem:onevar} now shows
that $(D,\sigma)$ is right $\sigma$-algebraically closed. For $n=1$
this is the whole conclusion.

Let $n\ge2$. Fix $a\in D^\times$ and put
$\theta=\Int(a^{-1})\circ\sigma$. By Lemma~\ref{lem:inner}, the
finite-dimensional module hypothesis is unchanged by this inner twist.
Moreover, $E_a(D,\sigma)=D^\theta$. Choose a pairwise commuting
tuple $b_1,\ldots,b_r\in D^\theta$ of length $r=n-2$, and set
\[
C=\Cen_{D^\theta}(b_1,\ldots,b_r).
\]
We must prove that this division ring is right algebraically closed.
Suppose otherwise. Lemma~\ref{lem:onevar}, with the identity
automorphism of $C$, gives a monic irreducible polynomial
\[
p=z^m+p_{m-1}z^{m-1}+\cdots+p_0\in C[z],
\qquad m\ge2.
\]

\smallskip\noindent
\emph{Step 1: Construct a finite-dimensional test module.}
Take the left $D[z]$-module
\[
M=D[z]/D[z]p.
\]
It has $D$-dimension $m$, and each coset has a unique representative
of degree less than $m$. Extend $\theta$ coefficientwise to $D[z]$
by fixing $z$, and define
\[
S(\bar q)=\overline{\theta(q)},\qquad
B_j(\bar q)=\overline{qb_j},\qquad
Z(\bar q)=\overline{zq},\quad\text{ for each }\bar{q}\in M.
\]
We verify these definitions and their relations before using them.
Every coefficient of $p$ lies in $D^\theta$, so $\theta(p)=p$ and
\[
\theta(hp)=\theta(h)p\in D[z]p.
\]
The same holds for $\theta^{-1}$. Thus $S$ is well-defined, bijective,
and $\theta$-semilinear. Each $b_j$ commutes with every coefficient
of $p$ and with the central variable $z$, so $pb_j=b_jp$. Hence
$(hp)b_j=hb_jp\in D[z]p$, proving that $B_j$ is well-defined.
The map $B_j$ is $D$-linear. The map $Z$ is the action of $z$ on
the quotient, and is $D$-linear because $z$ is central.

The maps $B_j$ commute pairwise because the $b_j$ do, and they
commute with $Z$. They also commute with $S$, since
$\theta(b_j)=b_j$; the identity $SZ=ZS$ follows from $\theta(z)=z$.
Consequently the following $r+2=n$ maps are commuting
$\theta$-semilinear operators:
\[
S,\ B_1S,\ldots,B_rS,\ ZS.
\]
For example, $(B_iS)(B_jS)=B_iB_jS^2=B_jB_iS^2$.
They define an $R_n(D,\theta)$-module structure on $M$.

\smallskip\noindent
\emph{Step 2: Recover coefficients in the smaller division ring.}
Choose a nonzero submodule of $M$ of  least $D$-dimension. It is simple,
so the hypothesis makes it a line $L=D\bar q$ for some $q\in D[z]$. Normalize a nonzero
reduced representative $q$ to be monic and write $k=\deg q<m$.
The map $S$ preserves $L$ and is injective. It therefore maps this
line onto itself, and $S^{-1}(L)=L$. Since each variable action
preserves $L$, the identities
\[
B_j=(B_jS)S^{-1},\qquad Z=(ZS)S^{-1}
\]
show that all $B_j$ and $Z$ preserve $L$.

Stability of $L$ under $S$ gives a scalar $d\in D$ with
$\theta(q)-dq\in D[z]p$. The difference has degree less than $m$,
so it is zero. Comparing the coefficient of $z^k$ in
$\theta(q)=dq$ gives $d=1$, because $q$ is monic. Thus every
coefficient of $q$ is fixed by $\theta$.
For each $j$, stability of $L$ under $B_j$ similarly gives an actual
polynomial identity $qb_j=c_jq$ for some $c_j\in D$. The coefficient
of $z^k$ gives $c_j=b_j$, and the remaining coefficients show that
all coefficients of $q$ commute with every $b_j$. We conclude that
\[
q\in C[z].
\]
When $r=0$, there are no $B_j$ to consider, and the fixed-coefficient
argument already gives this conclusion.

\smallskip\noindent
\emph{Step 3: Force a factorization over $C$.}
Because $Z$ preserves $L$, there is $c\in D$ such that
$(z-c)q\in D[z]p$. This is a nonzero monic polynomial of degree
$k+1\le m$. Every nonzero left multiple of $p$ has degree at least
$m$, so $k=m-1$. Its leading coefficient is $1$, which forces the
left multiplier of $p$ to be $1$. Therefore
\[
p=(z-c)q.
\]
Writing $q=z^{m-1}+q_{m-2}z^{m-2}+\cdots$ and comparing the
coefficient of $z^{m-1}$ yields
\[
c=q_{m-2}-p_{m-1}\in C.
\]
Both factors belong to $C[z]$ and have positive degree, contradicting
irreducibility of $p$. Hence every irreducible polynomial in $C[z]$
is linear, and Lemma~\ref{lem:onevar} gives right algebraic closure
of $C$.

The choices of $a$ and of the commuting tuple were arbitrary.
Together with the one-variable closure proved at the start, this is
exactly $\Hcond_n(D,\sigma)$.
\end{proof}
\section{Centrally finite algebraic closures}\label{sec:compression-tools}
\begin{proposition}\label{prop:finitefixed}
	Assume $\Hcond_2(D,\sigma)$. If every $E_a(D,\sigma)$ is finite-dimensional
	over its own center for every $a\in D^{\times}$, then $\Hcond_n(D,\sigma)$ holds for every $n\geq 1$.
\end{proposition}
\begin{proof}
	Fix $a\in D^\times$ and write $E=E_a(D,\sigma)$. It is right algebraically
	closed. If it is commutative, it is an algebraically closed field. Otherwise,
	the Niven--Jacobson--Baer classification of centrally finite right algebraically
	closed division rings says that $E$ is a Hamilton quaternion division algebra
	$\HH_R$ over a real closed field $R$; see \cite[Theorem 2.3]{Paran}.
	
	%In an algebraically closed field, every centralizer is the field itself. 
	
	Now let $\mathbf b=(b_1,\ldots,b_r)$ 	
	be any finite pairwise commuting tuple in \(\mathbb H_R\). If at least \(b_j\) is noncentral for some $j$, then 
	\[
	\mathrm C_{\mathbb H_R}(b_j)=F,\quad\text{where }F:=R(b_j).
	\]
	Moreover,
	\[
	{
		\mathrm C_{\mathbb H_R}(\mathbf b)=F=R(b_j).
	}
	\]
	 The field \(F=R(b_j)\) is  algebraically closed. 
	 
	 Otherwise, assume that every \(b_j\) is central. Then $\mathrm C_{\mathbb H_R}(\mathbf b)=\mathbb H_R$.
	  By the classification hypothesis, \(\mathbb H_R\) is right algebraically closed.  	 Thus $E$ is centrally algebraically closed; %compare 	\cite[Proposition 4.2]{AryCentral}.
	   This holds for every $a$, so all the conditions
	$\Hcond_n$ hold.
\end{proof}

\begin{lemma}\label{lem:fixeddimension}
	Let $F=Z(D)$ and assume $[D:F]<\infty$. For every automorphism $\theta$ of $D$,
	\[
	[D^\theta:F^\theta]\leq [D:F].
	\]
	In particular, $D^\theta$ is centrally finite.
\end{lemma}
\begin{proof}
	Suppose $e_1,\ldots,e_r\in D^\theta$ are linearly independent over $F^\theta$.
	They are also linearly independent over $F$. Indeed, otherwise choose an
	$F$-linear relation with minimal nonzero support, relabel its terms, and
	normalize one coefficient to obtain
	\[
	e_1+f_2e_2+\cdots+f_se_s=0,\qquad f_i\in F.
	\]
	Applying $\theta$ and subtracting gives a relation with strictly smaller
	support unless every $f_i$ is fixed by $\theta$. Minimality therefore implies
	$f_i\in F^\theta$ for all $i$, contrary to the assumed independence.
	It follows that $r\leq[D:F]$. Finally, $F^\theta\subseteq Z(D^\theta)$.
\end{proof}

\begin{lemma}[Weighted primitive element]\label{lem:weighted}
	Let $L/K$ be a finite separable field extension, where $K$ is infinite. Given
	$w\in L^\times$ and an integer $q\geq1$, there is $f\in L^\times$ such that
	\[
	L=K(wf^{-q}).
	\]
\end{lemma}
\begin{proof}
	By the primitive element theorem, write $L=K(t)$, for some $t\in L$; see
	\cite[§17.8, p. 208]{Stewart}. Let $\tau_1,\ldots,\tau_d$ be its distinct
	$K$-embeddings into an algebraic closure of $K$, where $d=[L:K]$.
	For $i\ne j$, the equality
	\[
	\frac{\tau_i(w)}{(\tau_i(t)-\lambda)^q}
	=\frac{\tau_j(w)}{(\tau_j(t)-\lambda)^q}
	\]
	excludes only finitely many $\lambda\in K$. In fact, after clearing
	denominators it is the vanishing of
	\[
	\tau_i(w)(\tau_j(t)-X)^q-
	\tau_j(w)(\tau_i(t)-X)^q.
	\]
	This polynomial is not identically zero: comparing leading coefficients
	would first give $\tau_i(w)=\tau_j(w)$, and then equality of the two
	powers would give $\tau_i(t)=\tau_j(t)$, a contradiction.
	Choose $\lambda\in K$ avoiding all these finitely many exceptions. Then $f=t-\lambda$ is nonzero and the $d$ images of
	$\alpha:=wf^{-q}$ are distinct.
	 Let $m_\alpha(X)\in K[X]$
	be the minimal polynomial of \(\alpha\) over \(K\). 	Since \(\tau_i\) fixes every coefficient of \(m_\alpha\),
	\[
	m_\alpha\bigl(\tau_i(\alpha)\bigr)
	=
	\tau_i\bigl(m_\alpha(\alpha)\bigr)
	=
	0.
	\]
	Therefore all \(d\) distinct elements
$\tau_1(\alpha),\ldots,\tau_d(\alpha)$
	are roots of \(m_\alpha\). Hence
	\[
	\deg m_\alpha\geq d.
		\]
	Thus $[K(wf^{-q}):K]=d$, or $K(wf^{-q})=L$.
\end{proof}

\begin{lemma}\label{lem:fixedimplementer}
	Suppose some positive power of $\sigma$ is inner. Then there exist $q\geq1$
	and $u\in D^\times$ with
	\[
	\sigma^q=\Int(u),\qquad \sigma(u)=u.
	\]
\end{lemma}
\begin{proof}
	Write $\sigma^m=\Int(v)$. We have
	\[
	\operatorname{Int}(\sigma(v))
	=
	\sigma\circ\operatorname{Int}(v)\circ\sigma^{-1}
	=\sigma\circ\sigma^m\circ\sigma^{-1}
	=\sigma^m
	=\operatorname{Int}(v).
	\]
	More generally, the same calculation gives, for every \(i\geq0\),
	\[
	{
		\operatorname{Int}(\sigma^i(v))
		=\sigma^i\circ\sigma^m\circ\sigma^{-i}
		=\sigma^m.
	}
	\]
	The equality
	$\Int(\sigma(v))=\Int(v)$ implies $\sigma(v)v^{-1}\in Z(D)$.
	Thus all elements $\sigma^i(v)$ are central multiples of $v$ and commute.
	Set
	\[
	u=v\sigma(v)\cdots\sigma^{m-1}(v).
	\]
	Since $\sigma^m(v)=v$, we have $\sigma(u)=u$. Moreover,
	$\Int(\sigma^i(v))=\sigma^m$ for every $i\ge 0$, so
	Therefore
	\[
	\operatorname{Int}(u)
	=
	\underbrace{\sigma^m\circ\cdots\circ\sigma^m}_{m\text{ factors}}
	=\sigma^{m^2}.
	\] Take $q=m^2$.
\end{proof}

\section{Proofs of main results}\label{sec:mainProofs}
\begin{proof}[{\bf Proof of Theorem~\ref{thm:main}}]
	Write $R=R_n(D,\sigma)$. We prove the module-theoretic equivalences
	first, then pass to ideals and zeros, and finally address the finite
	certificate condition.
	
	\smallskip\noindent
	\emph{Step 1: The  module-theoretic conditions.}
	Proposition~\ref{prop:sufficiency} gives
	\ref{main:H}$\Rightarrow$\ref{main:all}. Condition
	\ref{main:all} clearly implies \ref{main:finite}, since
	finite-dimensional simple modules form a subclass of all simple
	modules. Proposition~\ref{prop:necessity} gives
	\ref{main:finite}$\Rightarrow$\ref{main:H}. Thus these three
	conditions are equivalent: $(\text{i})\Leftrightarrow (\text{ii}) \Leftrightarrow (\text{iii}).$
	
	\smallskip\noindent
	\emph{Step 2: One-dimensional simple modules and point ideals.}
	Assume \ref{main:all}, and let $L$ be a maximal left ideal of $R$.
	Then $R/L$ is simple and one-dimensional over $D$. The nonzero vector
	$e=1+L$ is a $D$-basis, so there are unique $a_i\in D$ such that
	$x_ie=a_ie$. Commutativity of the variables gives
	\[
	\sigma(a_j)a_i e=x_i(x_je)=x_j(x_ie)=\sigma(a_i)a_j e.
	\]
	Independence of $e$ yields the admissibility equations, so
	$\ba=(a_1,\ldots,a_n)\in D^{n,\sigma}$. Also $(x_i-a_i)e=0$
	for every $i$, which means $\mathfrak m_{\ba}\subseteq L$.
	Proposition~\ref{prop:pointideal} makes $\mathfrak m_{\ba}$ a
	maximal proper left ideal. Hence $L=\mathfrak m_{\ba}$  We have proved \ref{main:ideals}.
	
	Now assume \ref{main:ideals} and let $M$ be any simple left
	$R$-module. A nonzero $v\in M$ generates it, so the map
	$R\to M$, $f\mapsto fv$, identifies $M$ with $R/L$ for some
	maximal left ideal $L$. By hypothesis $L=\mathfrak m_{\ba}$ for an
	admissible point, and Proposition~\ref{prop:pointideal} gives
	$R/L\simeq D_{\ba}$. Thus $\dim_DM=1$. This proves
	\ref{main:all}$\Leftrightarrow$\ref{main:ideals}.
	
	\smallskip\noindent
	\emph{Step 3: Maximal left ideals and common zeros.}
	Every proper left ideal in a unital ring is contained in a maximal
	left ideal. Under \ref{main:ideals}, a proper left ideal $I$
	therefore satisfies
	\[
	I\subseteq \mathfrak m_{\ba}
	\]
	for some admissible $\ba$. The kernel description in
	Proposition~\ref{prop:pointideal} says that every element of $I$
	vanishes at $\ba$. Hence \ref{main:zeros} holds.
	
	Conversely, assume \ref{main:zeros} and take a maximal left ideal
	$L$. A point $\ba\in V_\sigma(L)$ gives
	$L\subseteq\mathfrak m_{\ba}$. So $L=\mathfrak m_{\ba}$. Since every point ideal is already
	known to be maximal, this proves \ref{main:ideals}. We have proven $(\text{iv})\Leftrightarrow (\text{v}).$ Hence $(\text{i})\Leftrightarrow (\text{ii}) \Leftrightarrow (\text{iii})\Leftrightarrow (\text{iv})\Leftrightarrow (\text{v}).$
	
	\smallskip\noindent
	\emph{Step 4: Common zeros and finite certificates.}
	For a finite family $f_1,\ldots,f_r$, let
	$I=\sum_{j=1}^rRf_j$. The product formula
	\eqref{eq:product} gives
	\[
	V_\sigma(I)=V_\sigma(\{f_1,\ldots,f_r\}).
	\]
	Indeed, one inclusion follows because the generators belong to $I$;
	for the other, a common zero of the generators is a zero of every
	left multiple $g_jf_j$ and hence of every sum of such multiples.
	If this common zero set is empty and \ref{main:zeros} holds, $I$
	cannot be proper. Thus $1\in I$, which means
	$1=\sum_jg_jf_j$ for suitable $g_j\in R$. The reverse implication
	in \eqref{eq:certificate} holds without any closure assumption by
	Proposition~\ref{prop:pointideal}. This proves
	\ref{main:zeros}$\Rightarrow$\ref{main:certificate}.

	Assume \ref{main:certificate}, and let $I$ be proper. Since $R$ is Noetherian, we can choose a finite generating family as just established. If
	$V_\sigma(I)$ were empty, the equality of zero sets above and
	\ref{main:certificate} would yield $1\in I$, a contradiction.
	Thus \ref{main:zeros} follows. All six conditions are now equivalent.
	
	\medskip
	Now we prove the last statement. For any admissible point $\ba$, Proposition~\ref{prop:pointideal}
	gives the equivalence
	\[
	\ba\in V_\sigma(I)
	\quad\Longleftrightarrow\quad
	f(\ba)=0\text{ for every }f\in I
	\quad\Longleftrightarrow\quad
	I\subseteq\mathfrak m_{\ba}.
	\]
	Under $\Hcond_n$, it follows from Theorem~\ref{thm:main} that the family of maximal
	left ideals containing $I$ is exactly
	\[
	\{\mathfrak m_{\ba}:\ba\in V_\sigma(I)\}.
	\]
	Taking intersections and using the definition of $I_\sigma$ yields
	\[
	\bigcap_{\substack{L\supseteq I\\L\text{ maximal left}}}L
	=\bigcap_{\ba\in V_\sigma(I)}\mathfrak m_{\ba}
	=I_\sigma\bigl(V_\sigma(I)\bigr).
	\]
	
\end{proof}

\begin{proof}[{\bf Proof of Theorem~\ref{thm:embedding}}]
	\emph{Step 1: make the automorphism inner.}
	Start with the skew polynomial ring
	\[
	A_0=D[u;\sigma],\qquad ud=\sigma(d)u\quad(d\in D),
	\]
	where $u$ is an indeterminate. Since $\sigma$ is an automorphism,
	$A_0$ is a left and right Euclidean domain and has an Ore division
	ring of fractions; see \cite{Ore}. Denote this division ring by
	$K=D(u;\sigma)$. The natural map $A_0\to K$ is injective, so we
	regard $D$ and $u$ as lying in $K$. In particular, $u\ne0$ is
	invertible in $K$, and the defining relation gives
	\[
	udu^{-1}=\sigma(d)\qquad(d\in D).
	\]
	Thus conjugation by $u$ on $K$ already restricts to $\sigma$ on $D$.
	
	\smallskip
	\emph{Step 2: Pass to a centrally algebraically closed extension.}
	the embedding
	theorem of Aryapoor
	\cite[Theorem~5.2]{AryCentral}:
	a division ring $K$ has a centrally algebraically closed extension $E$
	for which $Z(E)=Z(K)$.
	
	Applying this theorem for the $K$ constructed in step 1 gives
	a centrally algebraically closed division ring $E\supseteq K$ such that
	\begin{equation}\label{eq:embedding-center}
		Z(E)=Z(K).
	\end{equation}
	Define $\widetilde{\sigma}(c)=ucu^{-1}$ for every $c\in E$.
	This is an automorphism of $E$, with inverse $c\mapsto u^{-1}cu$.
	The relation \(ud=\sigma(d)u\) gives \(\widetilde{\sigma}|_D=\sigma\). Moreover, for every \(r\geq 1\),
	\[
	\widetilde{\sigma}^{\,r}=\mathrm{id}_E
	\quad\Longleftrightarrow\quad
	u^r\in Z(E)=Z(K)
	\quad\Longleftrightarrow\quad
	\sigma^r=\mathrm{id}_D,
	\]
	since \(u^r\) commutes with \(u\) and \(K\) is generated as a division ring by \(D\) and \(u\). Thus the two automorphisms have the same order.
	
	\smallskip
	\emph{Step 3: Transfer skew variable to the central variable.}
	Consider $A=E[x;\widetilde{\sigma}]$ and put $y=u^{-1}x$ in $A$.
	For every $c\in E$,
	\[
	yc=u^{-1}xc
	=u^{-1}(ucu^{-1})x
	=cy.
	\]
	Hence $y$ commutes with all coefficients. The assignments
	$y\mapsto u^{-1}x$ and $x\mapsto uy$ define an isomorphism $$E[x,\widetilde{\sigma}];\longrightarrow E[y;\theta=\operatorname{int}(u^{-1})\circ\widetilde{\sigma}=\id_E]$$
	
	Since \(E\) is right algebraically closed, Lemma 2.5 implies that \((E,\widetilde{\sigma})\) is right \(\widetilde{\sigma}\)-algebraically closed.
	
	For \(a\in E^\times\), put \(v=u^{-1}a\). Then
	\[
	E_a(E,\widetilde{\sigma})
	=\{c\in E:ucu^{-1}a=ac\}
	=\mathrm C_E(v).
	\]
	For every finite pairwise commuting tuple \(\mathbf b\) in \(\mathrm C_E(v)\), including the empty tuple, \((v,\mathbf b)\) is a commuting tuple in \(E\), and
	\[
	\mathrm C_{\mathrm C_E(v)}(\mathbf b)
	=\mathrm C_E(v,\mathbf b)
	\]
	is right algebraically closed. Hence \(\mathrm C_E(v)\) is centrally algebraically closed, proving the remaining assertion.
\end{proof}

\begin{proof} [{\bf Proof of Corollary~\ref{cor:embedding-alln}}]
	Choose $E$, $u$, and $\widetilde{\sigma}$ as in
	Theorem~\ref{thm:embedding}. Its third assertion is exactly
	$\Hcond_1(E,\widetilde{\sigma})$.
	For $n\ge2$, fix $a\in E^\times$ and a pairwise commuting
	$(n-2)$-tuple in $E_a(E,\widetilde{\sigma})$.
	The fourth assertion says that this twisted fixed division ring is
	centrally algebraically closed. The centralizer of the chosen tuple is
	therefore right algebraically closed. Together with the third assertion,
	this proves $\Hcond_n(E,\widetilde{\sigma})$. Theorem~\ref{thm:main} now gives all the asserted
	Nullstellensatz conclusions. Innerness and preservation of order are
	part of the construction in Theorem~\ref{thm:embedding}.
\end{proof}

\begin{proof}[{\bf Proof of Corollary~\ref{cor:alln}}]
	Suppose first that the Nullstellensatz holds for every $n$.
	The case $n=1$ and Theorem~\ref{thm:main} give right
	$\sigma$-algebraic closure of $(D,\sigma)$. Fix $a\in D^\times$
	and a finite commuting tuple $\bb$ of length $r\ge0$ in $E_a$.
	Apply the theorem with $n=r+2$. Its centralizer condition says that
	$\Cen_{E_a}(\bb)$ is right algebraically closed. Since the tuple was
	arbitrary, including $r=0$, the division ring $E_a$ is centrally
	algebraically closed.
	
	Conversely, suppose $(D,\sigma)$ is right
	$\sigma$-algebraically closed and every $E_a$ is centrally
	algebraically closed. The first condition is $\Hcond_1$.
	For any fixed $n\ge2$, central algebraic closure of $E_a$ applies
	in particular to every commuting tuple of length $n-2$ in $E_a$.
	Thus $\Hcond_n$ holds for each $n$, and Theorem~\ref{thm:main}
	gives the  Nullstellensatz in every number of variables.
\end{proof}

\begin{proof}[{\bf Proof of Theorem~\ref{thm:compression}}] We prove the first part: 	
	
	For $a\in D^\times$, let $\theta_a=\Int(a^{-1})\circ\sigma$. The restrictions
	of $\theta_a$ and $\sigma$ to $F=Z(D)$ coincide. Lemma~\ref{lem:fixeddimension}
	gives
	\[
	[E_a(D,\sigma):F^\sigma]\leq[D:F],
	\]
	so every $E_a(D,\sigma)$ is centrally finite. Apply
	Proposition~\ref{prop:finitefixed}.
	
	\smallskip
For the second part, we first prove: for every $a\in
D^\times$ and every finite pairwise commuting tuple
$\mathbf b=(b_1,\ldots,b_r)$ in $E_a(D,\sigma)$, there is $f\in
E_a(D,\sigma)^\times$ such that
\[
\Cen_{E_a(D,\sigma)}(\mathbf b)=E_{af}(D,\sigma).
\]

Choose $q,u$ as in Lemma~\ref{lem:fixedimplementer} and put $K=F^\sigma$.
The restriction of $\sigma$ to $F$ has finite order. The fixed-field theorem
therefore gives $[F:K]<\infty$; see \cite[Theorem 10.5]{Stewart}. Thus $K$
is infinite. It is perfect: in positive characteristic, the unique $p$th
root in $F$ of any element of $K$ is also fixed by $\sigma$.
Since $D$ is algebraic over $F$ and $F/K$ is finite, every element of $D$
is algebraic over $K$.

Fix $a$ and $\mathbf b$, and set
\[
\theta=\Int(a^{-1})\circ\sigma,\qquad
w=a^{-1}\sigma(a^{-1})\cdots\sigma^{q-1}(a^{-1})u.
\]
Direct computation gives
\[
\theta^q=\Int(w),\qquad \theta(w)=w.
\]
Indeed, the second equality uses $\sigma(u)=u$ and
$\sigma^q(a^{-1})=ua^{-1}u^{-1}$. For every $x\in D^\theta$,
$\theta^q(x)=x$ implies $wx=xw$. Consequently,
\[
w\in Z(D^\theta)=Z(E_a(D,\sigma)).
\]
The field
\[
L=K(w,b_1,\ldots,b_r)
\]
is therefore commutative, lies in $E_a(D,\sigma)$, and is a finite separable
extension of $K$. Apply Lemma~\ref{lem:weighted} to obtain $f\in L^\times$
such that $L=K(d)$, where $d=wf^{-q}$.

Now put
\[
\psi=\Int(f^{-1})\circ\theta
=\Int((af)^{-1})\circ\sigma.
\]
Because $\theta(f)=f$ and $f$ commutes with $w$,
for every \(x\in D\),
\[
\theta\bigl(f^{-1}xf\bigr)
=\theta(f)^{-1}\theta(x)\theta(f)
=f^{-1}\theta(x)f.
\]
Thus
\[
\theta\circ\Int(f^{-1})
=\Int(f^{-1})\circ\theta=\psi.
\]
We can therefore take the \(q\)-th power of their composition:
\[
\begin{aligned}
	\psi^q
	&=\bigl(\Int(f^{-1})\circ\theta\bigr)^q\\
	&=\Int(f^{-1})^q\circ\theta^q\\
	&=\Int(f^{-q})\circ\Int(w)\\
	&=\Int(f^{-q}w)\\
	&=\Int(d).
\end{aligned}
\]

If $x\in D^\psi$, then $x$ commutes with $d$. Since $K\subseteq Z(D)$ and
$L=K(d)$, $x$ commutes with every element of $L$, in particular with $f$ and
all $b_i$. The equation $\psi(x)=x$ then gives $\theta(x)=x$, or $x\in E_a(D,\sigma)$. Thus
\[
D^\psi\subseteq\Cen_{D^\theta}(b_1,\ldots,b_r).
\]
Conversely, an element of $\Cen_{D^\theta}(b_1,\ldots,b_r)$ commutes with
$w$, because $w\in Z(D^\theta)$, and hence with all of $L$ and with $f$.
It is therefore fixed by $\psi$. This proves the asserted equality:
\[
\Cen_{E_a(D,\sigma)}(\mathbf b)=E_{af}(D,\sigma).
\]
Under $\Hcond_2(D,\sigma)$, the ring $E_{af}(D,\sigma)$ is right
algebraically closed. Hence every $E_a$ is centrally algebraically closed,
which proves the final assertion.
\end{proof}

\section{Specializations and examples}\label{sec:applications}

\subsection{Commutative coefficient fields}
\begin{corollary}\label{cor:field}
Let $k$ be a field and let $\sigma\in\Aut(k)$. For $n = 1$, the weak Nullstellensatz for $k[x; \sigma]$ holds if and
only if $(k, \sigma)$ is right $\sigma$-algebraically closed.   For any fixed
$n\ge 2$, the weak Nullstellensatz for $k[x_1,\ldots,x_n;\sigma]$
holds if and only if $(k,\sigma)$ is right
$\sigma$-algebraically closed and $k^\sigma$ is algebraically closed.
%If it holds for a value $n\ge2$, it holds for every $n\ge 2$.
\end{corollary}
\begin{proof}
	For $a\in k^\times$ and $c\in k$, commutativity and cancellation
	of $a$ give
	\[
	\sigma(c)a=ac
	\quad\Longleftrightarrow\quad \sigma(c)=c.
	\]
	Hence $E_a=k^\sigma$, independently of $a$. Since this is a field,
	every element commutes with every tuple, and the centralizer of any
	commuting tuple in $E_a$ is again $k^\sigma$ itself. For a field,
	right algebraic closure is the ordinary algebraic closure condition:
	every nonconstant polynomial has a root.
	
	Thus for every $n\ge2$, the condition $\Hcond_n(k,\sigma)$ is the
	same pair of requirements: right $\sigma$-algebraic closure of
	$(k,\sigma)$ and algebraic closure of $k^\sigma$.
	Theorem~\ref{thm:main} proves the asserted equivalence. If it holds
	for one such $n$, that pair of requirements holds and hence gives
	the result for every $n\ge2$. Its first requirement also gives the
	one-variable case by Lemma~\ref{lem:onevar}.
\end{proof}

\subsection{Central variables and inner automorphisms}
\begin{corollary}\label{cor:central}
	For central variables, the equivalent conditions of
	Theorem~\ref{thm:main} hold for $D[x_1,\ldots,x_n]$ if and only if
	the centralizer of every commuting $(n-1)$-tuple in $D$ is right
	algebraically closed. For $n=1$, this means that $D$ is right
	algebraically closed. For all $n$, the criterion is central algebraic
	closure of $D$.
\end{corollary}
\begin{proof}
	For $\sigma=\id$, the ring $R_n$ is the central polynomial ring
	$D[x_1,\ldots,x_n]$, and
	\[
	E_a=\{c\in D:ca=ac\}=\Cen_D(a).
	\]
	When $n=1$, the condition $\Hcond_1$ is exactly right algebraic
	closure of $D$, which is also the condition on the centralizer of
	the empty tuple.
	
	Now let $n\ge2$. If $b_1,\ldots,b_{n-2}$ commute pairwise in
	$E_a$, they also commute with $a$. Directly from the definitions,
	\[
	\Cen_{E_a}(b_1,\ldots,b_{n-2})
	=\Cen_D(a,b_1,\ldots,b_{n-2}).
	\]
	Suppose $\Hcond_n(D,\id)$ holds, and take any commuting
	$(n-1)$-tuple in $D$. If one entry is nonzero, reorder the tuple
	so that this entry is first and call it $a$. The displayed equality
	identifies its centralizer with one of those occurring in $\Hcond_n$,
	so it is right algebraically closed. If every entry is zero, its
	centralizer is $D$. This is right algebraically closed by the
	one-variable part of $\Hcond_n$; equivalently, the same centralizer
	is obtained from the tuple $(1,0,\ldots,0)$.
	
	Conversely, suppose the centralizer of every commuting
	$(n-1)$-tuple is right algebraically closed. Applying this assumption
	to the zero tuple first shows that $D$ is right algebraically closed.
	Applying it to $(a,b_1,\ldots,b_{n-2})$ for $a\ne0$ and using the
	displayed equality gives every other condition in $\Hcond_n$.
	Theorem~\ref{thm:main} proves the fixed-$n$ assertion.
	
	Finally, as $n\ge1$ varies, tuples of length $n-1$ exhaust all
	finite commuting tuples, including the empty tuple. Their centralizer
	condition is precisely central algebraic closure of $D$, proving the
	last assertion.
\end{proof}

This recovers \cite[Theorems~4.3 and 4.4]{AryCentral}.
More generally, if $\sigma=\Int(u)$, the change of variables
$y_i=u^{-1}x_i$ makes all variables central over $D$.
Thus the all-variable weak Nullstellensatz for an inner automorphism
is equivalent to central algebraic closure of $D$.

\subsection{Puiseux series and an infinite-order automorphism}\label{subsec:puiseux}
Let
\[
F=\bigcup_{m\ge1}\CC((s^{1/m})),
\]
with the usual compatible choices of fractional powers. Fix
$\lambda\in\RR_{>0}\setminus\{1\}$, and define
\[
\sigma_\lambda\!\left(\sum_r c_rs^r\right)
   =\sum_r c_r\lambda^r s^r.
\]
The skew Newton--Puiseux theorem of Paran and Vo \cite{ParanVo}
asserts that $(F,\sigma_\lambda)$ is right
$\sigma_\lambda$-algebraically closed. Comparing coefficients shows
that
\[
F^{\sigma_\lambda}=\CC,
\]
since $\lambda^r\ne1$ for every nonzero rational $r$.
Corollary~\ref{cor:field} therefore applies in every number of
variables.

This example does not satisfy the nonconstant-center hypothesis
mentioned in the introduction. Indeed, if
$p(t)=\sum_jf_jt^j\in F[t;\sigma_\lambda]$ commutes with $s$, then
comparison of the coefficient of $t^j$ gives
\[
f_j\lambda^j s=sf_j.
\]
Thus $f_j=0$ for $j\ge1$. A constant commutes with $t$ precisely
when it lies in $F^{\sigma_\lambda}$. Consequently,
\[
Z\bigl(F[t;\sigma_\lambda]\bigr)=\CC.
\]
Also $\sigma_\lambda$ has infinite order, since
$\sigma_\lambda^j(s)=\lambda^j s\ne s$ for $j\ge1$.

\subsection{Why one-variable closure is not sufficient}
\begin{example}\label{ex:frobenius}
Let $k=\overline{\FF_p}$ and let $\sigma(c)=c^p$.
For $a\in k$ and $j\ge1$,
\[
N_j(a)=a^{1+p+\cdots+p^{j-1}}.
\]
Therefore right evaluation of a positive-degree skew polynomial
produces a nonconstant ordinary polynomial in $a$. Since $k$ is
algebraically closed, $(k,\sigma)$ is right
$\sigma$-algebraically closed. However, $k^\sigma=\FF_p$ is not
algebraically closed. Hence the weak Nullstellensatz fails for every
$n\ge2$.

The failure can be seen directly in dimension two. Choose
$\beta\in\FF_{p^2}\setminus\FF_p$, and on $M=k^2$ define
\[
T_1(v,w)=(w^p,v^p),\qquad
E(v,w)=(\beta v,\beta^p w),\qquad T_2=ET_1.
\]
The map $E$ commutes with $T_1$, because $\beta^{p^2}=\beta$.
Thus $T_1,T_2$ define an $R_2(k,\sigma)$-module.
If a nonzero proper submodule existed, it would be a $k$-line.
It would be stable under $T_1^{-1}$ as well as $T_1$, and hence
under $E=T_2T_1^{-1}$. Since $\beta\ne\beta^p$, an $E$-invariant
line is a coordinate axis. But $T_1$ exchanges the two axes.
Therefore $M$ is simple of $k$-dimension two.
\end{example}

\begin{example}[Quaternionic coefficients]\label{ex:quaternion}
Take $D=\HH$ over $\RR$ and $\sigma=\Int(i)$, with the usual
quaternionic units $i,j,k$. The division ring $\HH$ is right
algebraically closed, and a centralizer of a commuting tuple is either
$\HH$ or a complex subfield. In particular, $\HH$ is centrally
algebraically closed; see \cite[Proposition~4.2]{AryCentral}.
For $a\ne0$,
\[
E_a=\Cen_{\HH}(i^{-1}a).
\]
Hence every $E_a$ is centrally algebraically closed, and the weak
Nullstellensatz holds for every $n$.
Nevertheless, $Z(\HH)^\sigma=\RR$ is not algebraically closed.
Thus algebraic closure of the fixed part of the center is not a
necessary condition in the division-ring setting.

The point $(k,2k)$ is admissible. Since $\sigma(k)=-k$,
\[
x_1^2(k,2k)=\sigma(k)k=1,
\]
although $k^2=-1$. Consequently $x_1^2-1$ vanishes at this skew
point. This illustrates the distinction between right skew evaluation
and ordinary substitution.
\end{example}

\end{document}